\documentclass[a4paper]{article}

\usepackage{graphicx}
\usepackage{latexsym}

\usepackage{amsmath,amssymb,amsfonts}
\usepackage{algorithm,algorithmic}
\usepackage{url}

\def\Underline{\setbox0\hbox\bgroup\let\\\endUnderline}
\def\endUnderline{\vphantom{y}\egroup\smash{\underline{\box0}}\\}
\def\|{\verb|}

\begin{document}

\title{Acceleration of Multiple Precision Matrix Multiplication using Ozaki scheme}




\author{Taiga Utsugiri\thanks{Shizuoka Institute of Science and Technology} \and Tomonori Kouya}

\maketitle

\begin{abstract}
Optimized multiple precision basic linear computation, especially matrix multiplication, is crucial for solving ill-conditioned problems. The recently proposed Ozaki scheme, which implements accurate matrix multiplication using existing optimized low precision matrix multiplication, is known to be useful for multiple precision as well. In this paper, we implement fixed precision multi-component-way matrix multiplication using Ozaki scheme and show that in some cases it is faster than existing optimized matrix multiplications. We also show that arbitrary precision matrix multiplication using Ozaki scheme is also faster than Strassen matrix multiplication up to a certain precision.
\end{abstract}



\section{Introduction}

The applications of floating-point arithmetic, which are not limited to scientific and technical calculations, are expanding in the current global Digital eXchange (DX), so it is essential to select the length of the mantissa in floating-point number according to the scale and purpose of the calculation and to use software optimized for the computer environment. It is desirable to be able to use the software comfortably both in a supercomputer environment specialized for large-scale calculations and in a consumer-oriented environment, which is the main axis of DX.

In most cases, the floating-point arithmetics natively supported by the hardware that performs the operations, such as CPUs and GPUs, are limited to IEEE754-1985 binary32 (24 bits in binary mantissa, about 7 decimal digits) and binary64 (53 bits in binary, about 15 decimal digits). The use of binary16 (11 bits, 3 digits) is expanding for deep learning applications. The floating-point operations with more than 64 mantissa digits, which are required for adverse problems or round-off-error-sensitive problems, rely on software implementations.

Currently, a free, high-performance floating-point arithmetic library that allows arbitrary mantissa settings, is almost ones based on the MPN (Multiple Precision Natural number) kernel of GNU MP\cite{gmp}, which is optimized for a variety of CPU architectures. MPFR\cite{mpfr} is also based on the MPN. In order to develop softwares with MPFR, major programming environments including C/C++, Python, and Julia have already provided class libraries and modules in a user-friendly form. For problems that require a mantissa with more than a few hundred digits, the use of MPFR should be considered first.

However, for many mildly adverse problems that do not require that many digits, fixed-precision floating-point operations in the neighborhood of IEEE754 binary128 (113 bits, 34 digits) can be used, which are faster than MPFR for the same number of digits. Among them, by combining multiple binary32 and binary64 and using the error-free transformation technique, double-double (DD, 106bits, 32 digits), triple-double (TD, 159bits, 49 digits), and quadruple-double (QD, 212bits, 64 digits) can be obtained. There have been many studies and implementations of multi-component fixed-precision computation with multiple precision floating-point operations, such as QD library\cite{qd} by Bailey's et.al. , GQD\cite{gqd} and others.

For vector and matrix operations that are fundamental to numerical computation, BLAS (Basic Linear Algebra Subprogram) \cite{blas}-based optimization libraries are routinely used, with Intel Math Kernel on CPUs and cuBLAS included in CUDA on GPUs. For linear calculations of binary64 or lower precision, these libraries can be used to perform fast calculations in consumer computing environments.

For linear calculations that require multiple precision, MPLAPACK/MPBLAS(\cite{mplapack}) by Nakata, which incorporates QD and MPFR as well as GCC-supported \_float128, is well-known. As of 2023, it has become the standard for multiple precision BLAS.

We have already demonstrated the performance of our optimized multi-component-way linear computation using AVX2 in consumer x86\_64 environment, aiming to implement a multi-precision basic linear computation faster than MPBLAS. In addition, to reduce the number of time-consuming multi-precision operations compared to binary64, we introduced the Strassen algorithm to matrix multiplication and succeeded in further speed-up together with AVX2. These implementations have also been found to contribute to the speedup of LU decomposition.

Recently, however, Ozaki scheme, which achieves high-precision matrix multiplication using binary32 and binary64 optimized matrix multiplications, has been proposed and actively studied on CPUs. On CPUs, there is an implementation of Ozaki scheme by Mukunoki et al. using \_float128\cite{mukunoki_binary128}, and on GPUs, there is an implementation of Ozaki scheme by Nanai et al. \cite{nanai_ozaki}.

Therefore, we have implemented Ozaki scheme on CPUs using the binary64 matrix multiplication (DGEMM) of Intel Math Kernel for DD(\cite{utsugiri_hpc2022_dd}), TD(\cite{utsugiri_hpc2022_td}), and QD(\cite{utsugiri_hpc2022_td}), and compare its usefulness with the existing AVX2-ized implementation of the Strassen algorithm (hereinafter abbreviated as ``Strassen + AVX2") in a consumer x86\_64 environment as shown in \tablename\ref{table:comp_env}. We also implemented the TS (Triple-Single, 72 bits, 18 digits) Ozaki scheme on GPUs and compared its performance with that of simple matrix multiplication using TS operations. Unfortunately, we were unable to confirm the usefulness of Ozaki scheme using cuBLAS on a consumer GPU, but on a CPU, the results showed that it outperformed our implementation of the AVX2-ized Strassen algorithm. The results also indicate the possibility of parallelization to achieve higher performance when the number of division in Ozaki scheme is increased.

\begin{table}
	\begin{center}
		\caption{Computational environment for benchmark tests}
		\label{table:comp_env}
		\begin{tabular}{c|c}\hline
			CPU	& Intel Core i7 11700 (8C/16T) \\
			Memory	& 32GB \\
			GPU	& NVIDIA GeForce RTX 3070 \\
			OS	& Ubuntu 18.04.5 LTS \\
			CUDA & 11.0 \\
			Intel One API & 2021.5.0 \\ 
			ICPC Compiler Option & -fp-model precise -O3 -qmkl\\ 
			QD, MPFR & 2.3.22, 4.1.0 \\ \hline
		\end{tabular}
	\end{center}
\end{table}

Fixed-precision (DD, TD, QD) linear computation on the CPU was implemented using QD and QD-based TD arithmetic C++ class libraries, and on the GPU using GQD\cite{gqd} and GQD-based TS and D+S arithmetic class libraries. Ozaki scheme is implemented using Intel Math Kernel's DGEMM on the CPU and cuBLAS's SGEMM/DGEMM on the GPU.

We also applied Ozaki scheme to arbitrary-precision matrix multiplication using MPFR, and experimentally confirmed that the performance improvement of Ozaki scheme over Strassen matrix multiplication up to a certain precision can be achieved.

The contents of this paper are as follows.

First, an overview of the position of this study in the optimization research of multi-precision linear computation is given. Next, the algorithm of Ozaki scheme for multiple precision matrix multiplication is presented, and the results of the following benchmark tests are discussed.

\begin{enumerate}
	\item DD, TD, and QD Matrix multiplication and its application to LU decomposition
	\item D+S, DD, and, TS matrix multiplication on GPU
	\item MPFR arbitrary precision matrix multiplication on CPU
\end{enumerate}

We focus on the problem of finding the matrix product $C := AB$ for given real square matrices $A$ and $B$$\in\mathbb{R}^{n\times n}$. The elements of $A$ and $B$ used for benchmarking are pseudo-random numbers generated by (\ref{eqn:err1})

\begin{equation}
a_{ij}, b_{ij} := (ru - 0.5) \times \exp(rn), \label{eqn:err1}
\end{equation}
where $ru$ is a uniform random number of $[0, 1]$ and $rn$ is a random number following a standard normal distribution. The computed matrix product error represents the maximum relative error in all elements of computed $C$.
\section{Positioning of our study in optimized multiple precision linear computation}

As mentioned above, multiple precision floating point arithmetic is used when standard hardware-oriented floating-point operations such as binary16, binary32, and binary64, cannot maintain precision due to missing mantissa, or when there is a need to derive a high-precision formulas. 

Bailery et al.'s QD library \cite{qd} is well known as an implementation based on a multi-component method that uses multiple binary32 and binary64 numbers and maintains the precision of operations by using error-free transformation techniques. There are also libraries for multiple precision linear computation based on this library that have been optimized using parallelization techniques such as SIMD instructions, OpenMP, and MPI, including MPLAPACK/MPBLAS \cite{mplapack} and Lis \cite{lis}. There is also research on matrix multiplication using Ozaki scheme, such as Mukunoki et al.\cite{mukunoki_binary128} and Nanai et al.\cite{nanai_ozaki}, which are known to be faster than existing implementations. A list of current research on optimization of these multiple precision matrix multiplications is shown in \tablename\ \ref{table:opt_mp_blas}.

\begin{table}[htb]
\begin{center}
	\caption{List of Optimization Methods and Research for Multiple Precision Basic Linear Computation}
	\label{table:opt_mp_blas}
	\includegraphics[width=.75\textwidth]{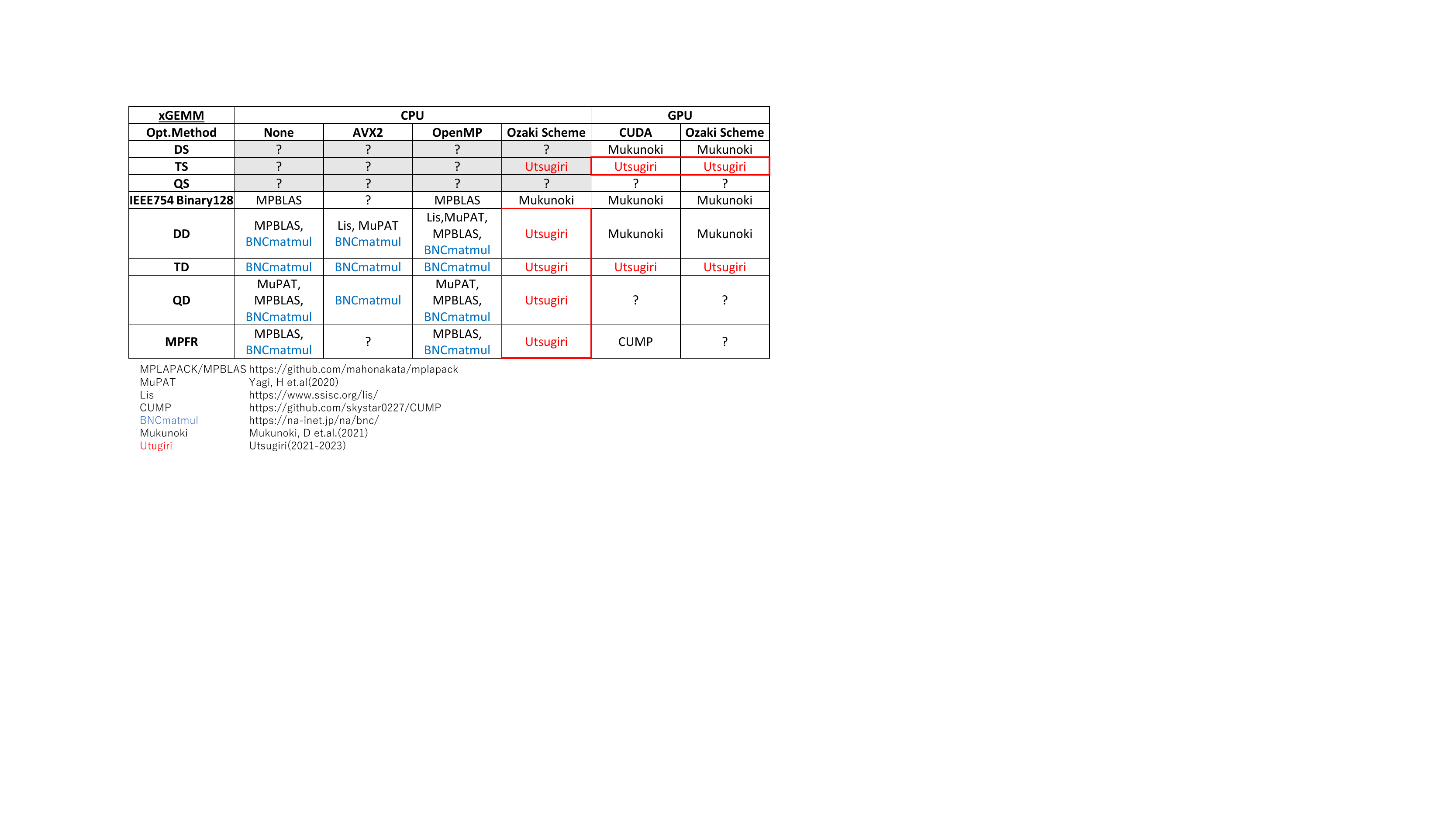}
\end{center}
\end{table}

There are no existing libraries for multiple precision linear computation that support optimization using Ozaki scheme and that are widely available free of charge. Considering the future expansion of the use of multiple precision numerical computation, the emergence of high-performance basic linear computation libraries optimized for hardware architectures, such as OpenBLAS\cite{openBLAS} and Intel Math Kernel, with fewer language and software environment restrictions, is required. We are focusing on the development of a high-performance basic linear computation library that is optimized for hardware architectures.

As shown in the red thick-bordered box in \tablename\ \ref{table:opt_mp_blas}, we focus on DD-, TD-, QD-, and MPFR-type matrix multiplications with Ozaki scheme on CPUs and TS-type matrix multiplications with Ozaki scheme on GPUs.

%
\section{Implementation of matrix multiplication using Ozaki scheme}

Ozaki scheme for multiple precision matrix multiplication is an algorithm that aims at both speed-up and accuracy improvement by dividing an original long-precision matrix into short-precision matrices, much like ``Split" function in error-free transformation technique, in order to take advantage of the speed-up of optimized short precision matrix multiplication (xGEMM) functions. For a given matrix $A \in \mathbb{R}^{m\times l}$, $B \in \mathbb{R}^{l\times n}$, to obtain a matrix product $C := AB \in \mathbb{R}^{m\times n}$ of long $L$-bit precision, $A$ and $B$ are divided into matrices of short $S$-bit precision ($S << L$), is shown in Algorithm\ \ref{algo:ozaki_scheme}. The $S$-bit algorithm is used for the operations not described above, and the $L$-bit algorithm is used only where high-precision operations are required.

Hereafter, the sets of $S$- and $L$-bit length binary floating-point numbers are denoted $\mathbb{F}_{bS}$ and $\mathbb{F}_{bL}$, respectively. Following this, the set of binary32 and binary64 floating-point numbers can be written as $\mathbb{F}_{b24}$ and $\mathbb{F}_{b53}$.

\begin{algorithm}[htb]
    \algsetup{linenosize=\small}
    \small
    \caption{Ozaki scheme for multiple precision matrix multiplication}\label{algo:ozaki_scheme}
    \hspace*{\algorithmicindent} \textbf{Input:} $A \in \mathbb{F}_{bL}^{m\times l}, B \in\mathbb{F}_{bL}^{l\times n}$ \\
	\hspace*{\algorithmicindent} \textbf{Output:} $C \in \mathbb{F}_{bL}^{m\times n}$
    \begin{algorithmic}
        \STATE $A^{(S)} := A$, $B^{(S)} := B$ : $A^{(S)} \in \mathbb{F}_{bS}^{m\times l}$, $B^{(S)} \in \mathbb{F}_{bS}^{l\times n}$
		\STATE $\mathbf{e} := [1\ 1\ ...\ 1]^T\in \mathbb{F}_{bS}^l$
		\STATE $\alpha := 1$
		\WHILE{$\alpha < D$}
			\STATE ${\boldsymbol\mu}_A := [\max_{1 \leq p \leq l} |A^{(S)}_{ip}|]_{i=1, 2, ..., m} \in \mathbb{F}_{bS}^m$
			\STATE ${\boldsymbol\mu}_B := [\max_{1 \leq q \leq l} |B^{(S)}_{qj}|]_{j=1, 2, ..., n} \in \mathbb{F}_{bS}^n$
			\STATE ${\boldsymbol\tau}_A := [2^{\lceil \log_2(({\boldsymbol\mu}_A)_i)\rceil + \lceil (S + \log_2(l)) / 2 \rceil} ]_{i = 1, 2, ..., m} \in \mathbb{F}_{bS}^m$
			\STATE ${\boldsymbol\tau}_B := [2^{\lceil \log_2(({\boldsymbol\mu}_B)_j)\rceil + \lceil (S + \log_2(l)) / 2 \rceil} ]_{j = 1, 2, ..., n} \in \mathbb{F}_{bS}^n$
			\STATE $S_A := \boldsymbol\tau_A \mathbf{e}^T$
			\STATE $S_B := \mathbf{e} \boldsymbol\tau_B^T$
			\STATE $A_\alpha := (A^{(S)} + S_A) - S_A$: $A_\alpha \in \mathbb{F}_{bS}^{m\times l}$
			\STATE $B_\alpha := (B^{(S)} + S_B) - S_B$: $B_\alpha \in \mathbb{F}_{bS}^{l\times n}$
			\STATE $A := A - A_\alpha$, $B := B - B_\alpha$ : $L$-bit computation
			\STATE $A^{(S)} := A$, $B^{(S)} := B$
			\STATE $\alpha := \alpha + 1$  
		\ENDWHILE
		\STATE $A_D := A^{(S)}$, $B_D := B^{(S)}$
		\STATE $C := O$
		\FOR{$\alpha = 1, 2, ..., D$}
			\FOR{$\beta = 1, 2, ..., D - \alpha + 1$}
				\STATE $C_{\alpha\beta} := A_\alpha B_\beta$
			\ENDFOR
			\STATE $C := C + \sum^{D - \alpha + 1}_{\beta = 1} C_{\alpha\beta}$ : $L$-bit computation
		\ENDFOR
    \end{algorithmic}
\end{algorithm}

\begin{figure}[htb]
	\begin{center}
		\includegraphics[width=.75\textwidth]{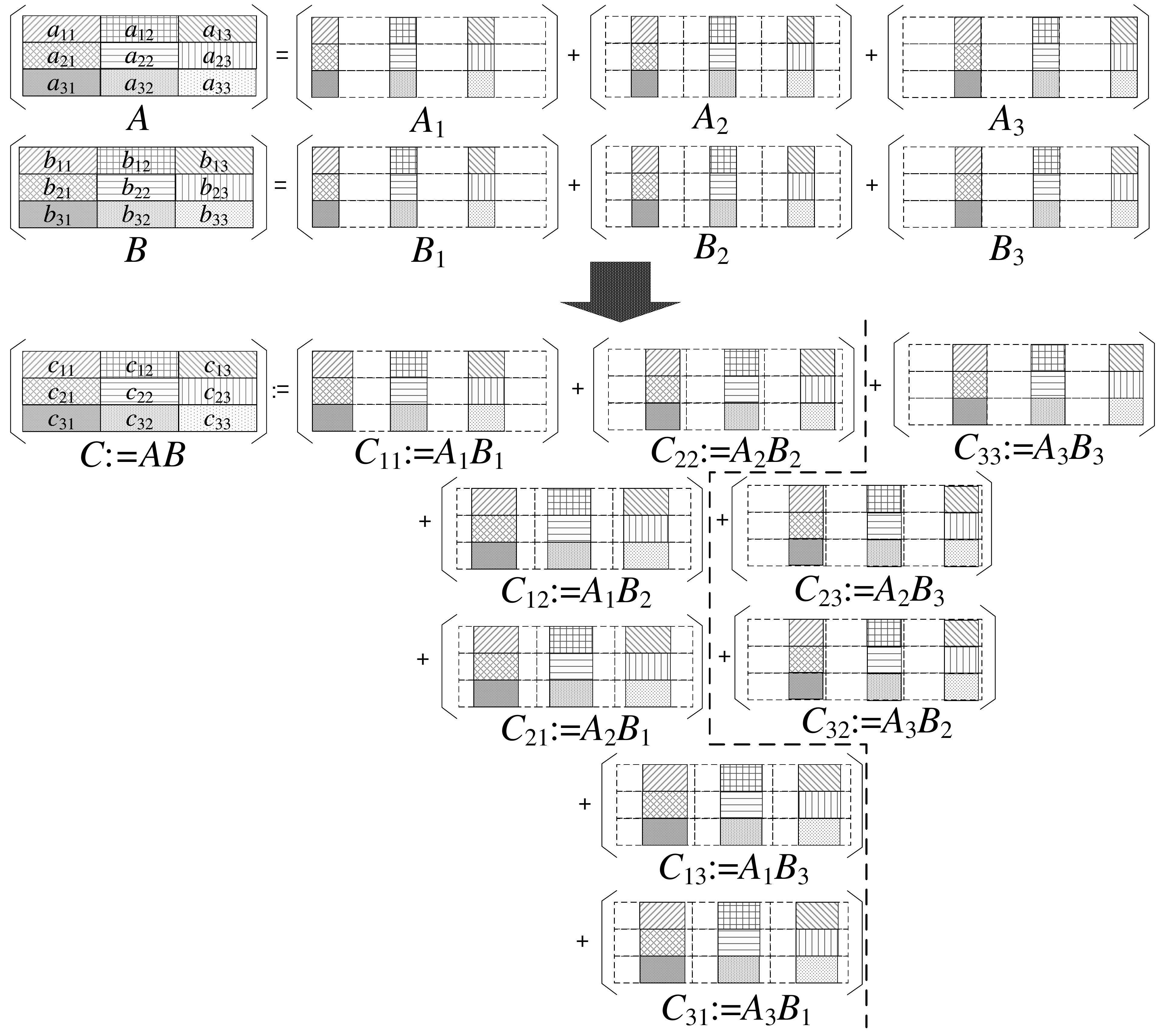}
		\caption{Matrix multiplication based on Ozaki scheme when divided into 3 parts}
		\label{fig:ozaki_scheme3x3}
	\end{center}
\end{figure}

A schematic representation of matrix multiplication when $A$ and $B\in\mathbb{R}^{3\times 3}$ are divided into 3 parts, respectively, is shown in\figurename\ \ref{fig:ozaki_scheme3x3}.

As a concrete example of \figurename\ \ref{fig:ozaki_scheme3x3}, consider the case where Ozaki scheme is executed using TD-type matrix multiplication ($L = 159$) with binary64 DGEMM ($S = 53$). In this case, each element of $A_\alpha$ and $B_\beta$ $(\alpha, \beta = 1, 2, 3)$ generated by partitioning $A$ and $B$ is divided into $C_{\alpha\beta} = A_\alpha B_\beta$ so that no error occurs in any $C_{\alpha\beta} = A_\alpha B_\beta$ calculation. The length of each elements of $A_\alpha$ and $B_\beta$, is determined to fit in the hypothetical director. Therefore, if the addition in $C = \sum_{\alpha,\beta} C_{\alpha\beta}$ is performed in TD type, each error-free term $C_{\alpha\beta}$ can be used to obtain $C\approx AB$ in TD precision. In practice, it is possible to remove more of the lower $C_{\alpha\beta}$, but in our benchmark test, all the obtained $C_{\alpha\beta}$ are added. As shown in \figurename\ \ref{fig:ozaki_scheme3x3}, $C_{11}$, $C_{12}$, $C_{13}$, $C_{21}$, $C_{22}$ and $C_{31}$ are obtained and added.

In summary, the pros and conds of Ozaki scheme are depending on
\begin{enumerate}
	\item Actual number of $A$ and $B$ partitions required,
	\item Short precision ($S$-bit) xGEMM performance used,
	\item Amount of loss of significant digits in $AB$ calculation.
\end{enumerate}
The results of the profile described below show that the number of digits in the $AB$ calculation is determined by the number of divisions used. According to the profile results described below, in the case of multiple precision, the performance of addition and subtraction is somewhat related when the matrix size is small and $L$ is large, but its dependence decreases as the matrix size increases. Therefore, these three factors have a significant impact on the performance of obtaining a high-precision matrix product.

In other words, if the distribution of matrix elements is short, the number of divisions is small, and if $S$-bit precision xGEMM is fast, it is expected to be faster than ordinary matrix multiplication, and fewer terms are needed if there are fewer losses of digits in the $AB$ calculation. Conversely, if the distribution of matrix elements is long, the number of divisions increases and the number of terms required to maintain accuracy also increases, and if there is a large number of loss of digits in the $AB$ calculation, more terms are required to maintain the accuracy of $C$. Thus, unless the properties of $A$ and $B$ are known in advance, "the performance cannot be known without running the real calculation," which is a characteristic of Ozaki's scheme.

Therefore, the benchmark test results shown below are only for the test matrix generated by the (\ref{eqn:err1}) formula, and should be considered to be subject to significant changes depending on $A$, $B$ and the computer environment. Nevertheless, previous studies have shown that there is much room to take advantage of fast xGEMM features such as Intel Math Kernel, cuBLAS, and MAGMA.

%
\section{Performance evaluation of fixed precision matrix multiplication}

As mentioned above, previous studies by Mukuraki et al. have shown that Ozaki scheme contributes to speeding up multiple precision matrix multiplication in the case of matrix multiplication without extreme loss of digits, but so far no comprehensive results have been presented for multi-component fixed-precision DD-, TD-, and QD-type matrix multiplications. We therefore evaluated the performance of Ozaki scheme on all currently available fixed-precision operations, including TS-type operations on GPUs, with a focus on triple-word operations, which were the subject of our earlier work, and compared them with relative errors for comparative purposes. The results are summarized here.

%
\subsection{Relative error and computational time of DD-type matrix multiplication}

\begin{figure}[htb]
	\begin{center}
		\includegraphics[width=.75\textwidth]{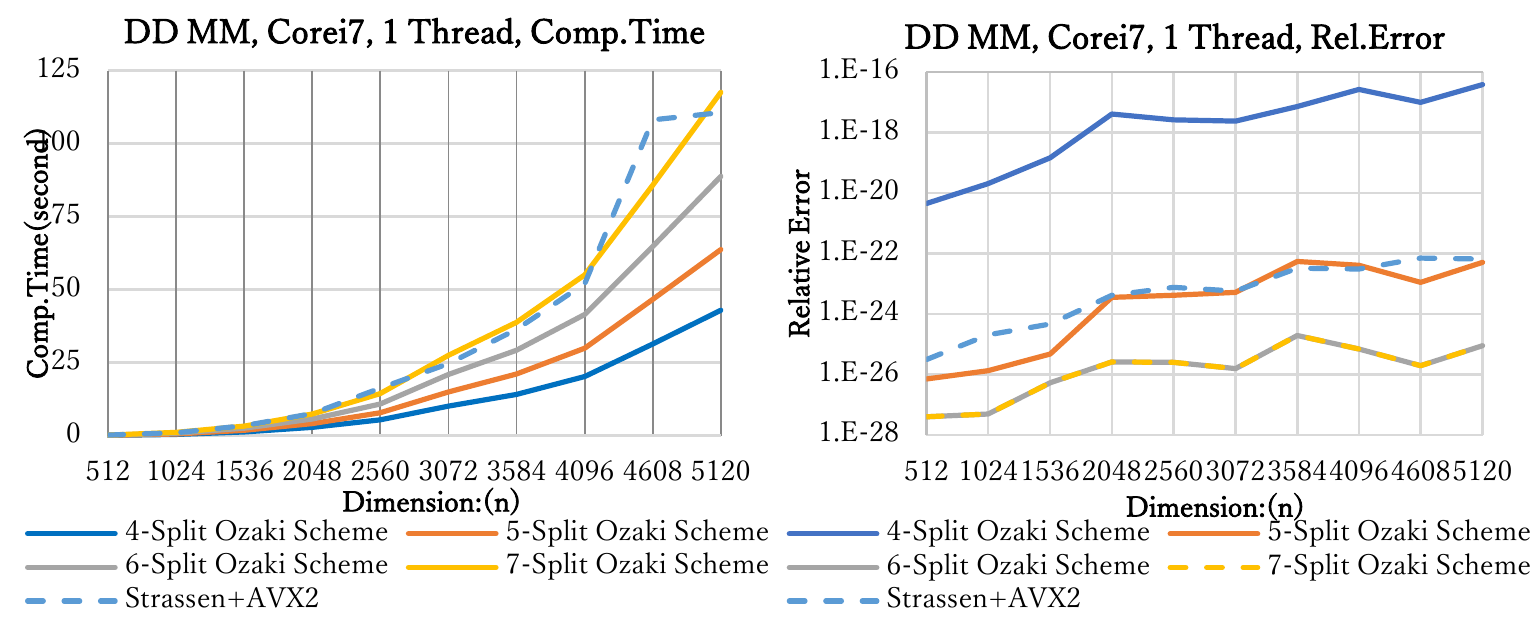}
		\caption{Computational time (left) and relative error (right) of DD matrix multiplication: $n=512, 1024, ..., 5120$}
		\label{fig:dd_1th_time_relerr}
	\end{center}
\end{figure}

The benchmark results for the DD Ozaki scheme on CPU are shown in \figurename\ \ref{fig:dd_1th_time_relerr}. In this figure, ``simple matrix multiplication" means using a triple loop to get a matrix product.

In \figurename\ \ref{fig:dd_1th_time_relerr}, the left figure shows the execution time and the right figure shows the relative error. The computation time for simple matrix multiplication is omitted from these figures because it is too large. Let us start with the right figure. The vertical axis shows the relative error and the horizontal axis shows the execution time. The number of dimensions run is from $n=512$ to $5120$. The figure shows that Ozaki scheme with 5 segments is as accurate as Strassen+AVX2. The maximum accuracy of Ozaki scheme is obtained after 6 divisions, and the accuracy does not improve as the number of divisions is increased.

The figure on the left is illustrated below. The vertical axis represents the execution time, and the horizontal axis represents the number of dimensions. This figure shows that the speeds of Ozaki scheme with 7 divisions and Strassen+AVX2 are equivalent. The 6-segmentation scheme, which maximizes the accuracy of Ozaki scheme, is about 1.2 times faster, and the 5-segmentation scheme, which achieves the same accuracy as Strassen+AVX2, is about 1.8 times faster.

\begin{figure}[htb]
	\begin{center}
		\includegraphics[width=.75\textwidth]{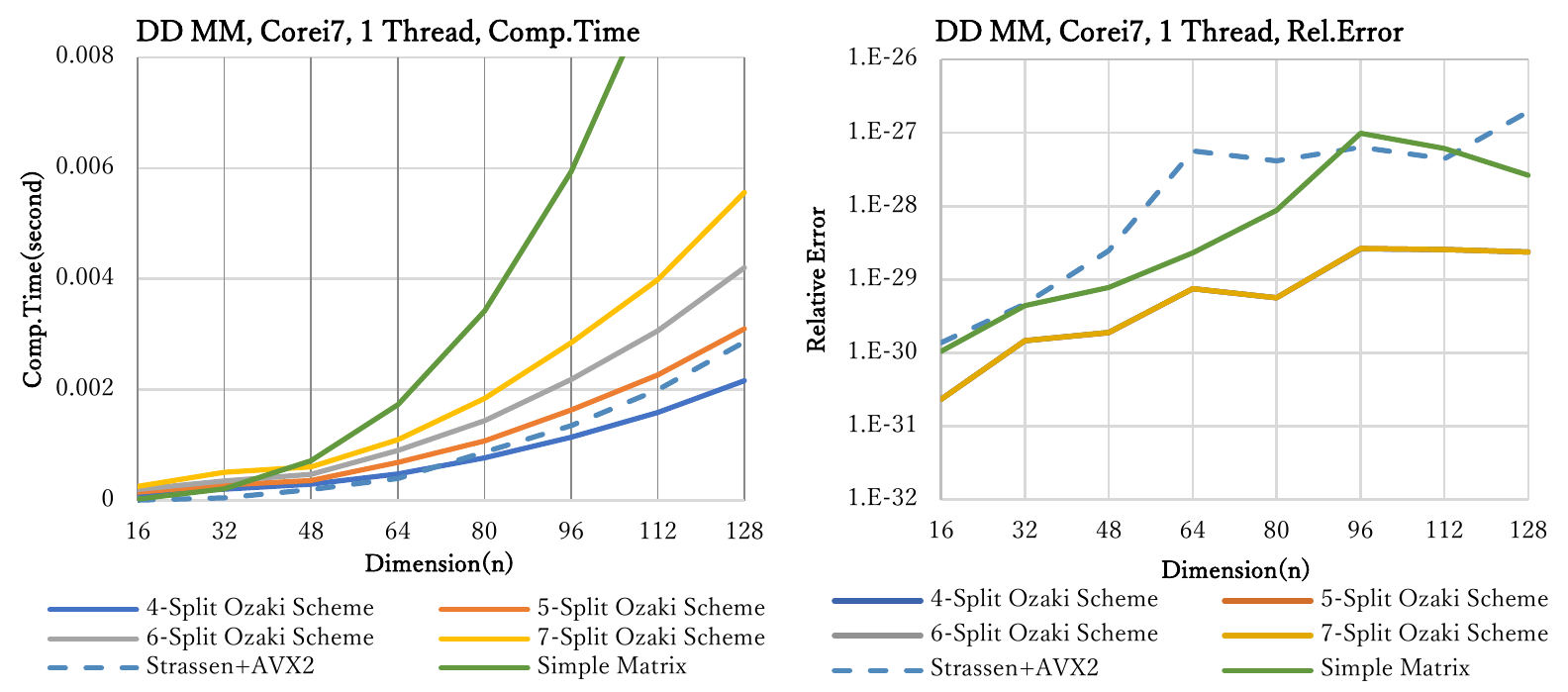}
		\caption{Computational time (left) and relative error (right) of DD matrix multiplication: $n=16, 32, ..., 128$}
		\label{fig:dd_1th_time_relerr_small_dim}
	\end{center}
\end{figure}

\figurename\ \ref{fig:dd_1th_time_relerr_small_dim} shows the benchmarks for dimension $n=16$ to $128$. For comparison, here again we compare the runtimes and relative errors for simple matrix multiplication and Strassen+AVX2.

First, let's start with the relative error shown on the right. The vertical axis represents the relative error, and the horizontal axis represents the number of dimensions. This figure shows that the accuracy is better than simple matrix multiplication and Strassen+AVX2 up to 128 dimensions for Ozaki scheme with 4 or more divisions.

Next, we explain the runtime shown in the left figure. The vertical axis represents the execution time, and the horizontal axis represents the number of dimensions. The figure shows that simple matrix multiplication is faster than Ozaki scheme up to 32 dimensions, and that Ozaki scheme is not faster than Strassen+AVX2 up to 128 dimensions, which is equivalent to Ozaki scheme of 5 divitions.

The profiling results of the above Ozaki scheme runs, showing the ratio of execution time for division, multiplication, and addition for the 5-division and 7-division cases, are shown in Figure 1. The left figure shows Ozaki scheme with 5 divisions, which is equivalent to Strasse+AVX2 in accuracy, and the right figure shows 7 divisions, which is the maximum number of divisions performed in this study.

\begin{figure}[htb]
	\begin{center}
		\includegraphics[width=.75\textwidth]{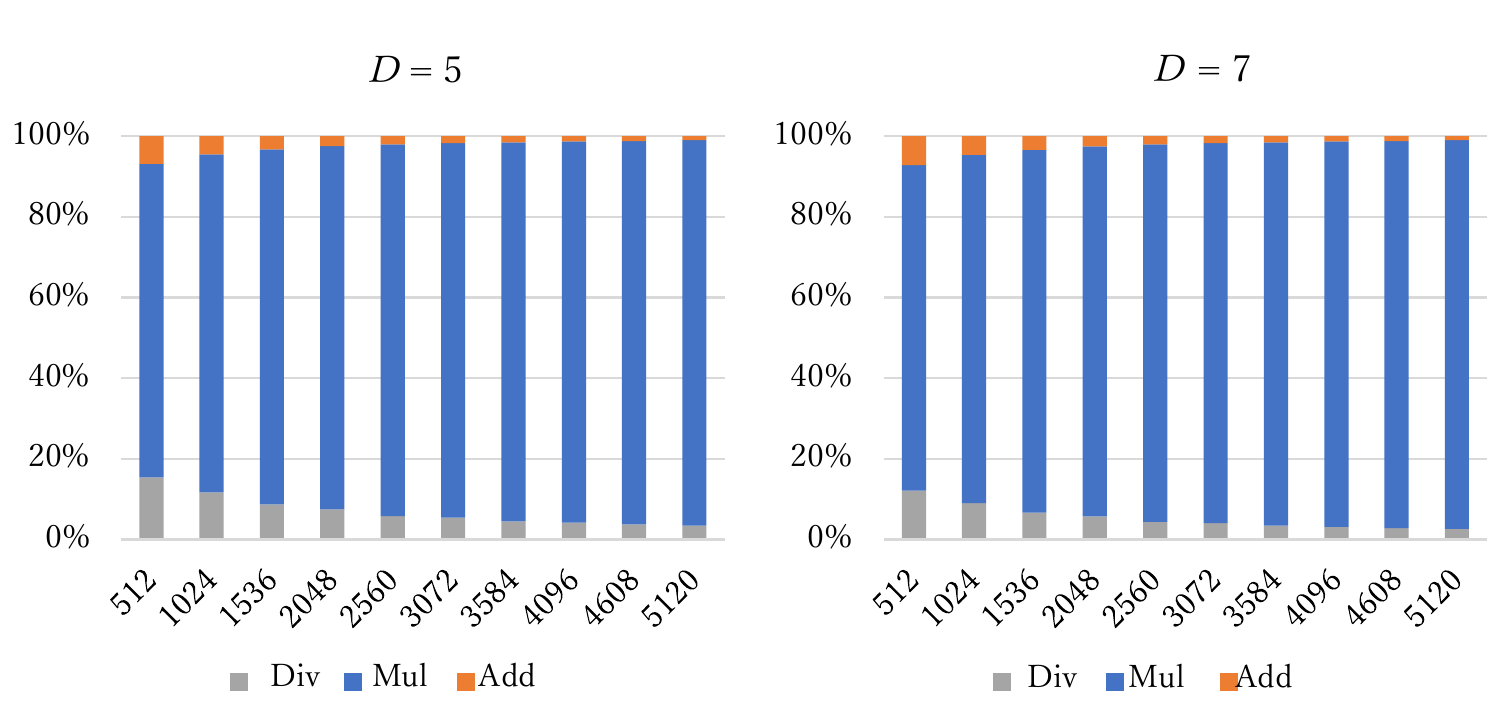}
		\caption{Profiling result of DD matrix multiplication: $D=5$(left) and $D=7$(right)}
		\label{fig:dd_profiling}
	\end{center}
\end{figure}

For both 5-division and 7-division schemes, it can be seen that as the matrix size increases, the operation time of the multiplication part increases more than the operation time of the division part and the addition part. Compared to the TD scheme (discussed below), the proportion of the multiplication part increases, and in the DD Ozaki scheme, most of the execution time is spent on multiplication. In case of more dimensions than 3072, the results show that matrix multiplication (DGEMM) shares more than 90\% of the operations in Ozaki scheme.

%
\subsection{Relative error and computational time of TD-type matrix multiplication}

\begin{figure}[htb]
\begin{center}
	\includegraphics[width=.75\textwidth]{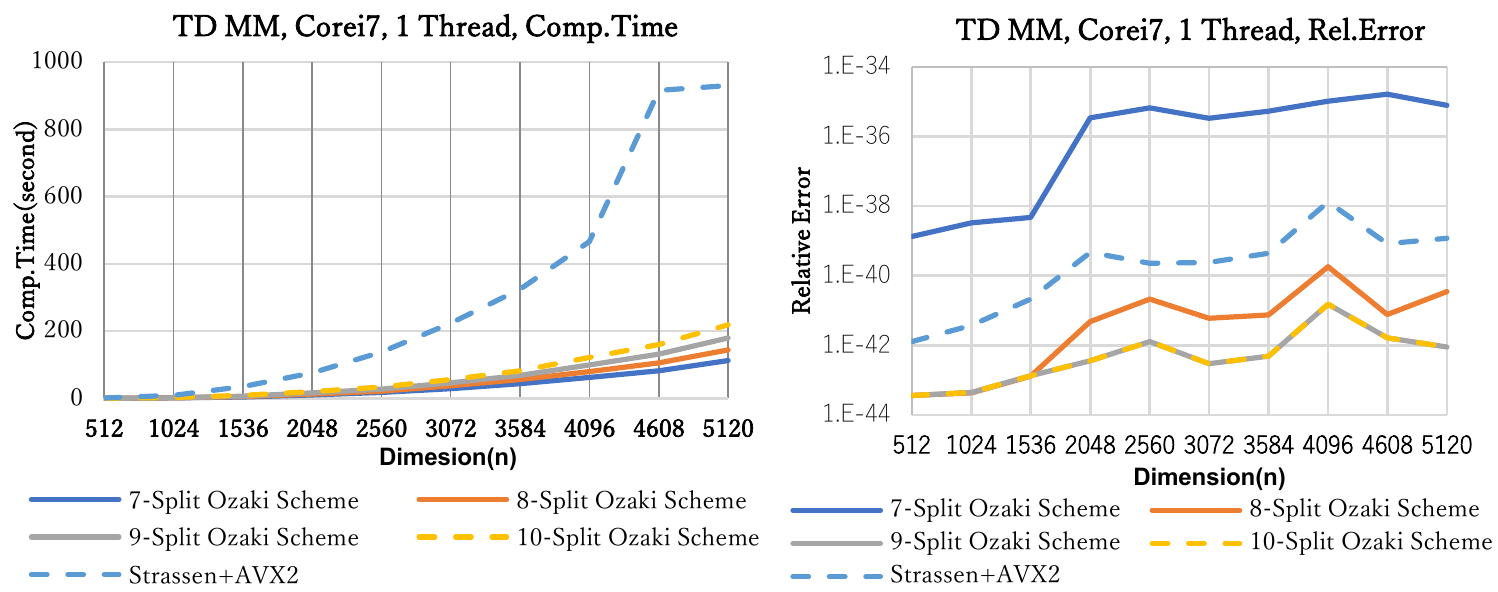}
	\caption{Computational Time (left) and relative error (right) of TD matrix multiplication: $n=512, 1024, ..., 5120$}
	\label{fig:td_time_relerr}
\end{center}
\end{figure}

Next, we show the results of benchmarking the execution time of matrix multiplication using the TD Ozaki scheme in \figurename\ \ref{fig:td_time_relerr}. The results have already been published in a previous paper, but Ozaki scheme is a new benchmark in a new consumer computer environment with fewer operations than in the previous paper.

First, let's start with the right figure of \figurename\ \ref{fig:td_time_relerr}. The right figure shows the relative error, where the vertical axis indicates the relative error and the horizontal axis indicates the number of dimensions. The number of dimensions ranges from $n=512$ to $5120$, with calculations performed every 512. The results of simple matrix multiplication are not shown in this figure because they are too large computational time.

The graph shows that TD Ozaki scheme is about 0.3 orders of magnitude less accurate at 9 divisions, and that the relative error is no smaller at 10 or more divisions. If you want the same or better accuracy than Strassen+AVX2, you need more than 8 divisions of Ozaki scheme.

Next, the left figure is explained. This figure shows the execution time, with the vertical axis representing the execution time and the horizontal axis representing the number of dimensions executed. The graph shows that the TD-type Ozaki scheme is approximately 5.7 times faster than Strassen+AVX2 by a factor of at most 5.7 for the 10 partitions where the accuracy of the TD-type Ozaki scheme is maximized, and approximately 8.6 times faster by a factor of at most 8 for Ozaki scheme with the same or higher accuracy than Strassen+AVX2.

\begin{figure}[htb]
	\begin{center}
		\includegraphics[width=.75\textwidth]{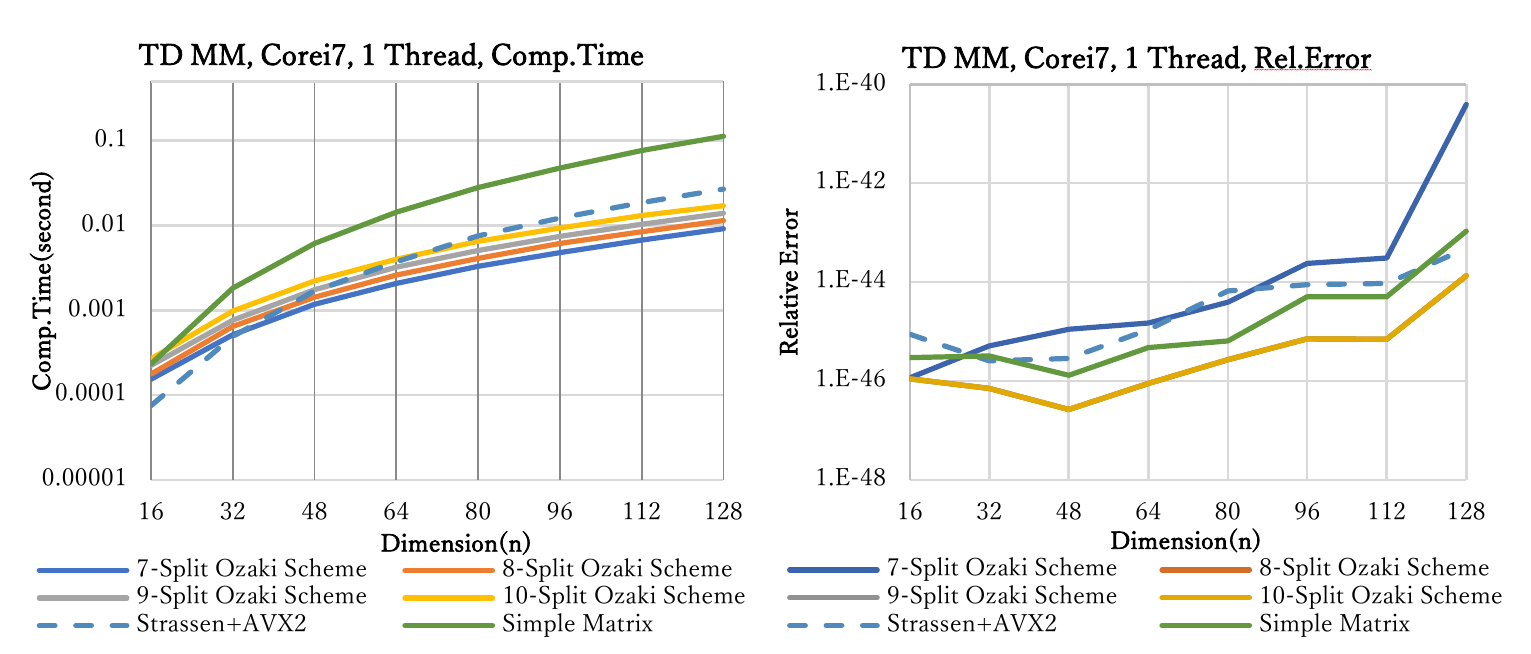}
		\caption{Computational Time (left) and relative error (right) of TD matrix multiplication: $n=16, 32, ..., 128$}
		\label{fig:td_time_relerr_small_dim}
	\end{center}
\end{figure}

\figurename\ \ref{fig:td_time_relerr_small_dim} shows the benchmark results with $n=16$ to $128$ . This compares the execution times and relative errors for Ozaki scheme, simple matrix multiplication, and matrix multiplication using Strassen+AVX2.

We begin with an explanation of the relative error diagram on the right-hand side. This figure shows that the relative errors of simple matrix multiplication and Strassen + AVX2 are included between the 7-division Ozaki scheme and the 8-division Ozaki scheme and above. It was found that Ozaki scheme gives better accuracy with a smaller number of divisions when the matrix size is smaller than in the case of the DD scheme shown earlier. This indicates that the number of divisions may rise with increasing dimension

Next, the left-hand side of the runtime diagram is explained. This figure shows that simple matrix multiplication is faster than Ozaki scheme for up to 16 dimensions, and Strassen+AVX2 is faster than Ozaki scheme for up to 64 dimensions. After that, Ozaki scheme was the fastest over simple matrix multiplication and Strassen+AVX2. This indicates that simple matrix multiplication and Strassen+AVX2 are faster than Ozaki scheme for small matrix sizes.

\begin{figure}[htb]
	\begin{center}
		\includegraphics[width=.75\textwidth]{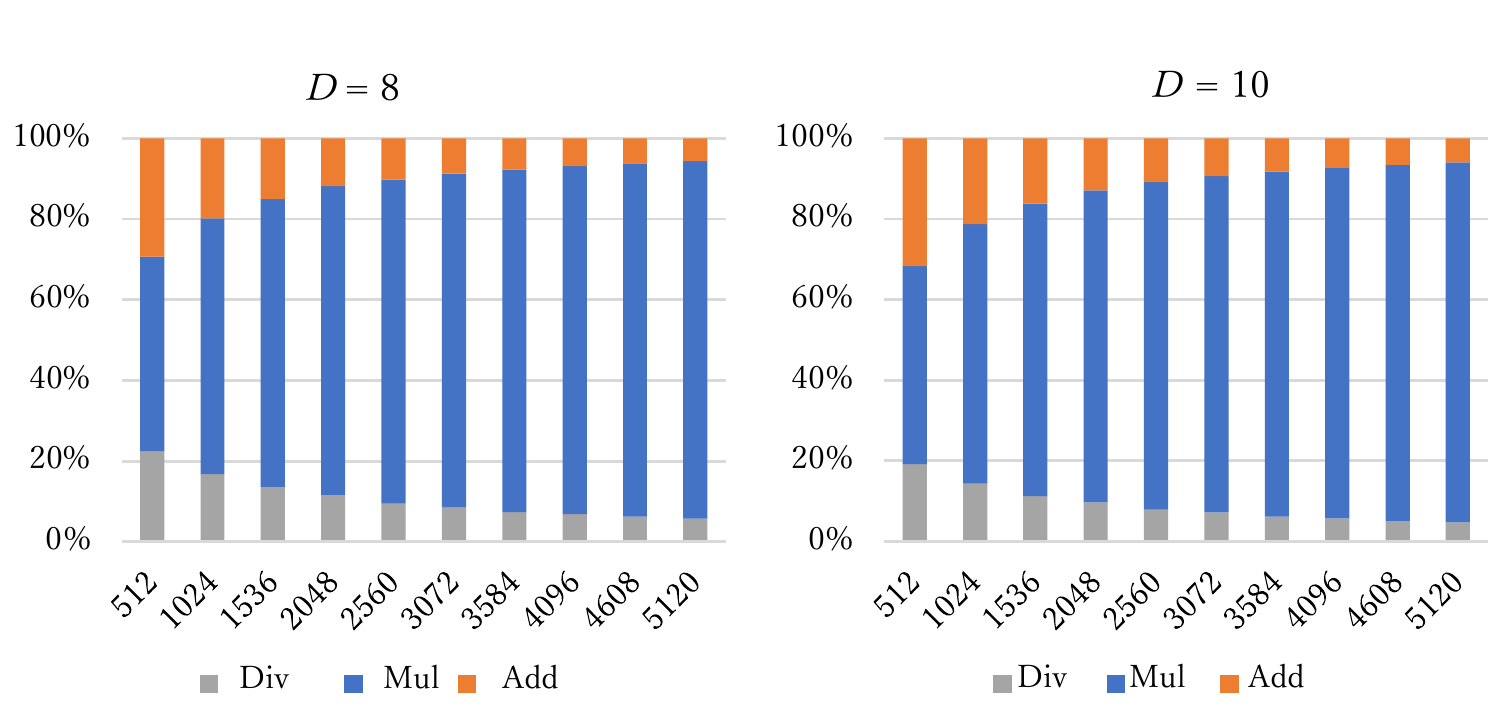}
		\caption{Profiling results of TD matrix multiplication: $D=8$(left) and $D=10$(right)}
		\label{fig:td_profiling}
	\end{center}
\end{figure}

\figurename\ \ref{fig:td_profiling} shows the runtime profiling results for the TD-type Ozaki scheme with 8 and 10 divisions.

The figure on the left shows an 8-division scheme that achieves accuracy better than Strassen+AVX2, and the figure on the right shows a 10-division scheme that maximizes the accuracy of Ozaki scheme at present. The right figure shows the 10-division scheme, which is currently the most accurate with respect to the current accuracy of Ozaki's scheme.

%
\subsection{Relative error and computational time of QD-type matrix multiplication}

\begin{figure}[htb]
	\begin{center}
		\includegraphics[width=.75\textwidth]{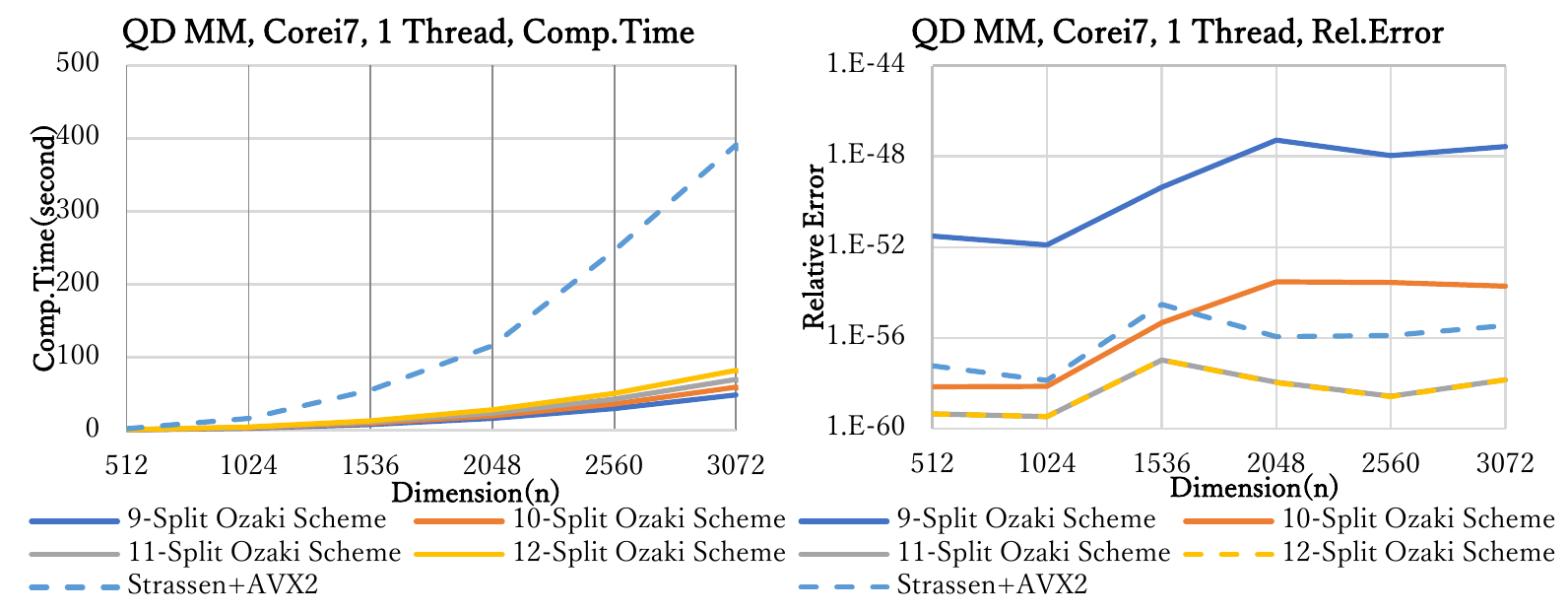}
		\caption{Computational time (left) and relative error (right) of QD matrix multiplication: $n=512, 1024, ..., 5120$}
		\label{fig:qd_1th_time_relerr}
	\end{center}
\end{figure}

\figurename\ \ref{fig:qd_1th_time_relerr} shows the performance evaluation of matrix multiplication using the QD-type Ozaki scheme on a CPU. The left side shows the execution time and the right side shows the relative error. For comparison, Strassen+AVX2 and simple matrix multiplication were used. Again, the simple matrix multiplication is omitted from the figure because its runtime is too large.

First, let us start with the relative error shown in the figure on the right. The vertical axis shows the relative error and the horizontal axis shows the number of dimensions, which were measured in 512-dimensional intervals from $n=512$ to $3072$. The figure shows that the accuracy of Ozaki scheme with 10 divisions and Strassen+AVX2 is equivalent in the range of 512 to 1536 dimensions. Later, when the number of dimensions increases to $2048$ or more, the accuracy of Strassen+AVX2 falls between the 10 and 11 divisions of Ozaki scheme. The number of divisions at which Ozaki scheme achieves maximum accuracy is found to be the 11 divisions of Ozaki scheme for dimensions up to 2560.

Next, we discuss the execution time in the left figure. The vertical axis shows the execution time and the horizontal axis shows the number of dimensions. This figure shows that Ozaki scheme is about 6.8 times faster than Strassen+AVX2 for 10 divisions, and about 5.7 times faster for 11 divisions.

\begin{figure}[htb]
	\begin{center}
		\includegraphics[width=.75\textwidth]{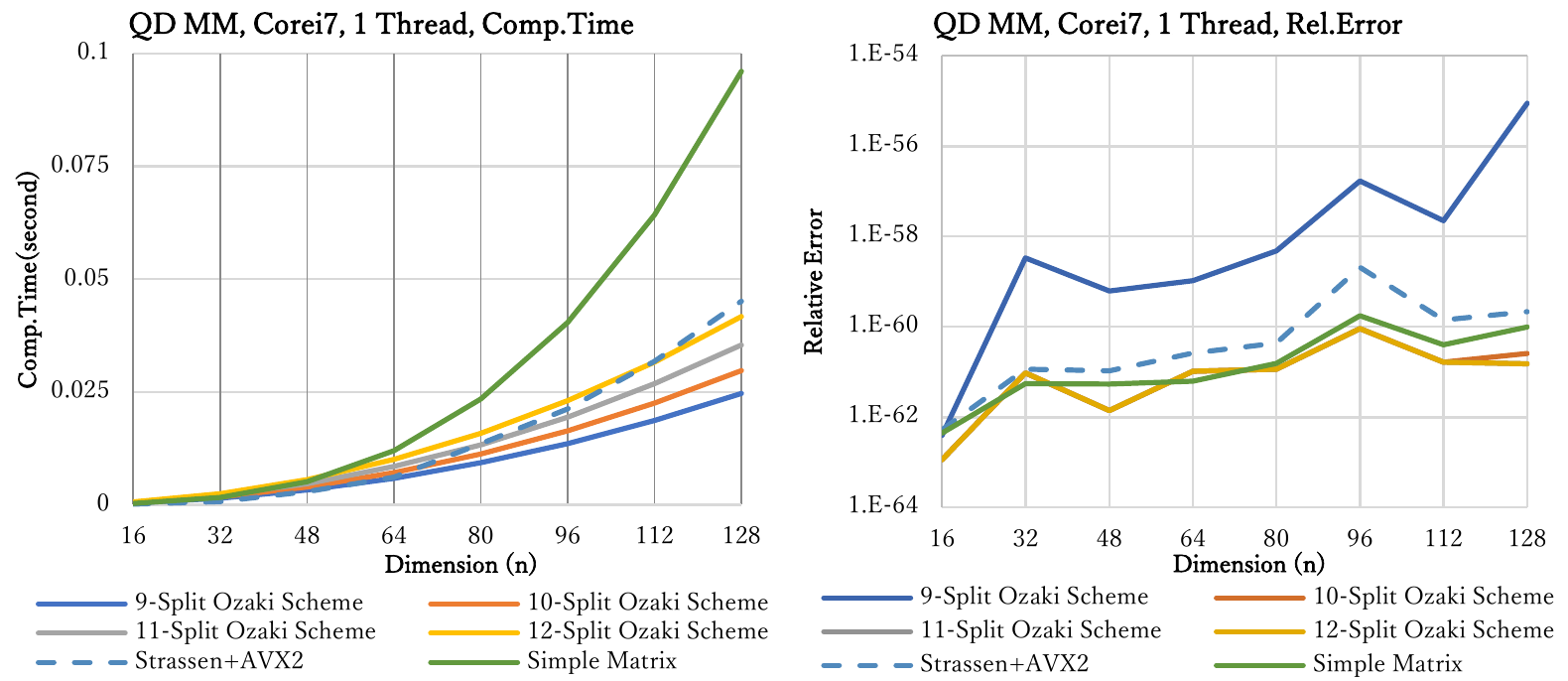}
		\caption{Computational time (left) and relative error (right) of QD matrix multiplication: $n=16, 32, ..., 128$}
		\label{fig:qd_1th_time_relerr_small_dim}
	\end{center}
\end{figure}

\figurename\ \ref{fig:qd_1th_time_relerr_small_dim} shows the benchmark results for dimension $n=16$ to $128$. For comparison, we use the execution time and relative error of simple matrix multiplication and Strassen+AVX2.

First, let's start with the figure on the right. The vertical axis shows the relative error and the horizontal axis shows the number of dimensions. This figure shows that the accuracy of simple matrix multiplication and Strassen+AVX2 are between the 9- and 10-segment accuracy of Ozaki scheme. The left graph shows the execution time, where the vertical axis indicates the execution time and the horizontal axis indicates the number of dimensions. The figure shows that simple matrix multiplication is faster than 12-division Ozaki scheme up to 48 dimensions, and faster than Strassen+AVX2 up to 112 dimensions.

%
\subsection{TS-type matrix multiplication on GPU}

Nanai et al.'s results demonstrate the usefulness of Ozaki scheme on consumer GPUs. How is the performance of Ozaki scheme for higher-precision TS-type operations (triple-single, 72 bits, and 21 digits)? And how much does the TS-type operation on GPUs improve on the D+S operation of Mukuraki et al.? To answer these questions, we present the results of our benchmark tests below.

\begin{figure}[htb]
	\begin{center}
		\includegraphics[width=.35\textwidth]{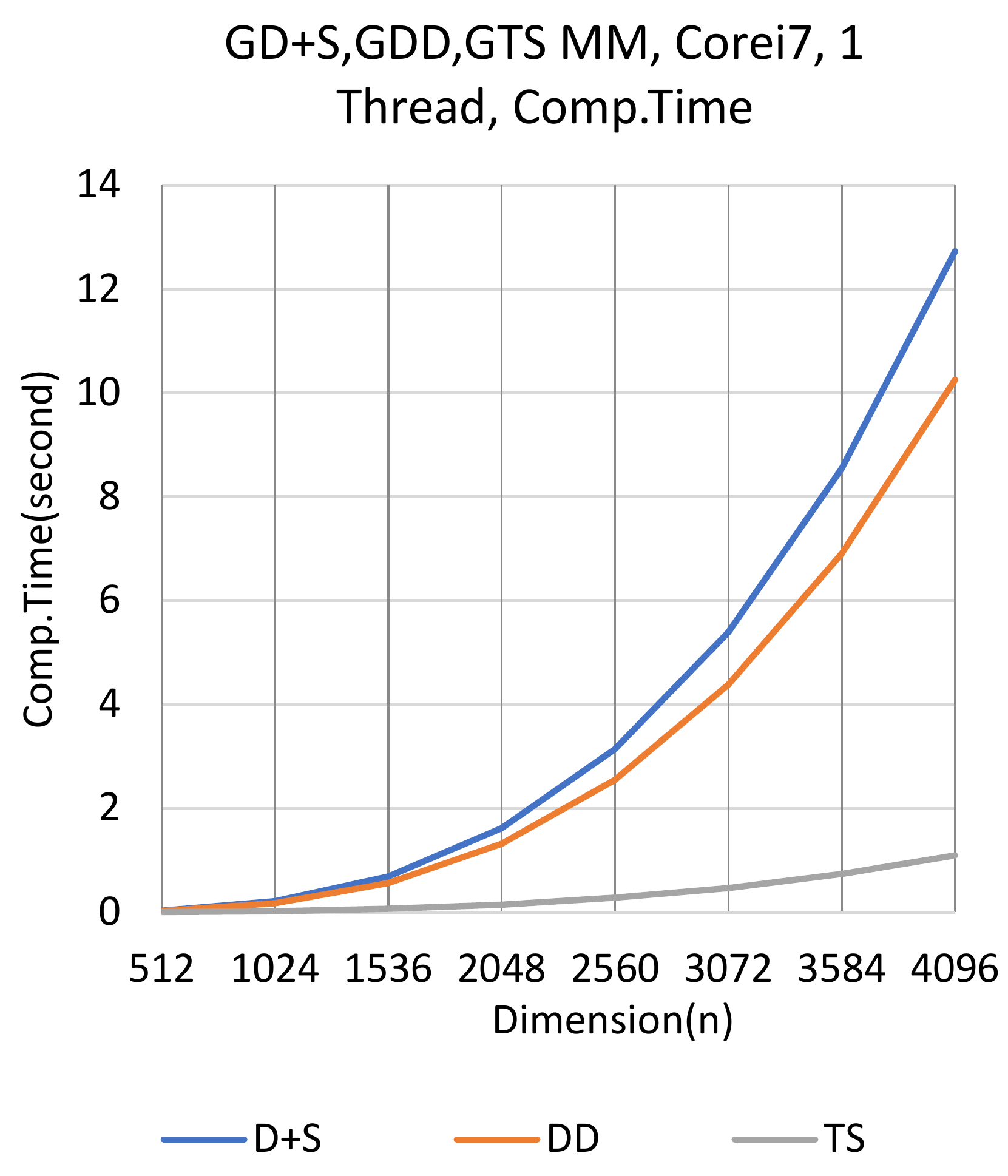}
		\includegraphics[width=.35\textwidth]{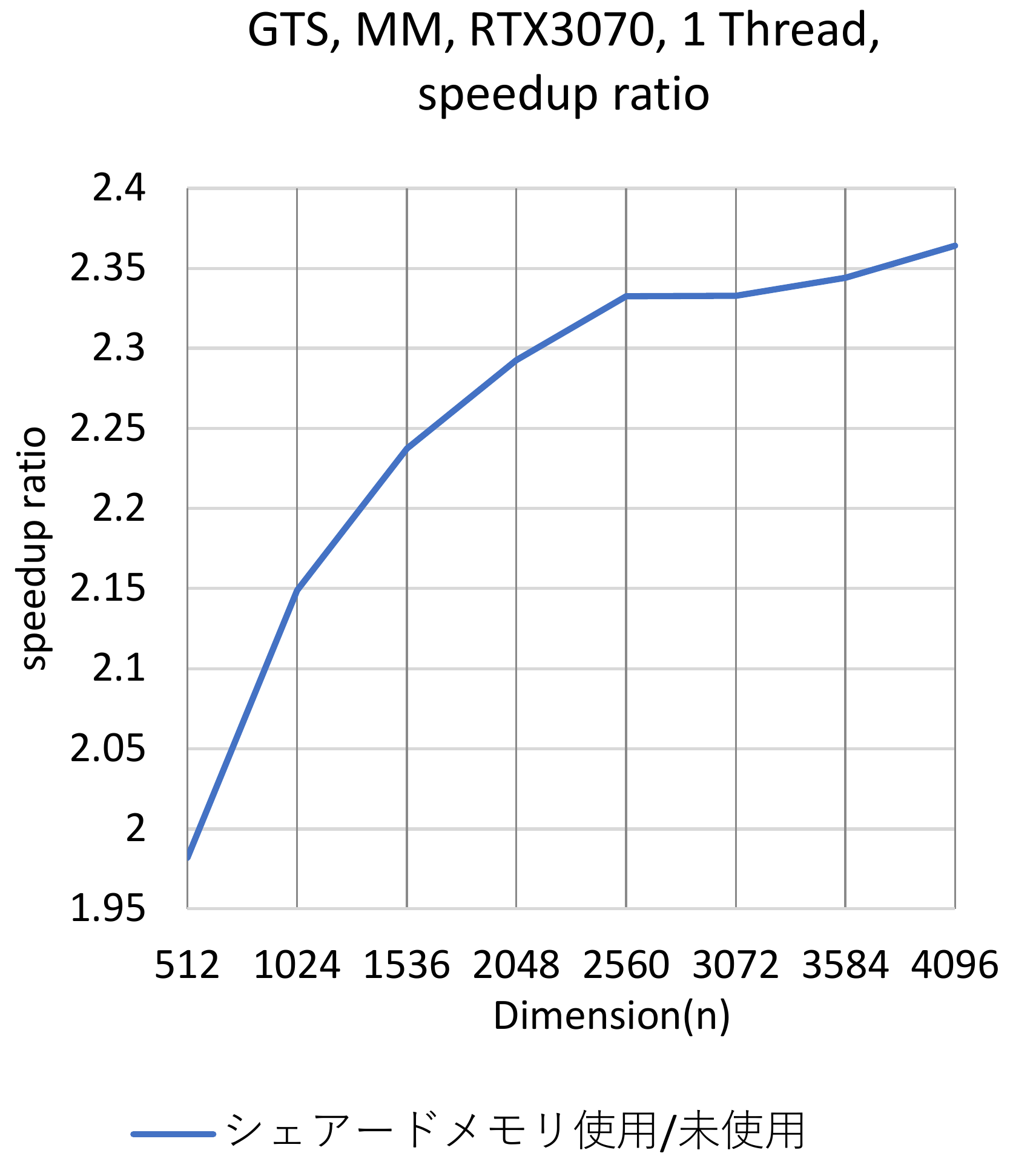}
		\caption{Computational time of D+S, DD, TS matrix multiplication(left) and speedup ratio using shared memory (right)}
		\label{fig:gpu_ds_dd_ts_time}
	\end{center}
\end{figure}

The left figure of \figurename\ \ref{fig:gpu_ds_dd_ts_time} shows the performance evaluation results of TS-type matrix multiplication on a GPU. The difference from our previous work is that we compared DD-type matrix multiplication with D+S-type matrix multiplication and improved Ozaki scheme. Matrix sizes ranged from $n=512$ to $4096$, measured in increments of 512. The vertical axis shows the running time and the horizontal axis shows the number of dimensions.

The graph shows that the TS type is faster than the DD and D+S types. The reason for the similar speeds of the DD and D+S types is that the D+S and DD types have virtually the same algorithm. Furthermore, D+S is slightly slower than DD because it converts the lower-level binary32(single prec.) number to binary64 (double prec.) type.

Next, a comparison of execution time multipliers with and without GPU shared memory is shown in the right figure of \figurename\ \ref{fig:gpu_ds_dd_ts_time}.

The vertical axis shows the execution time multiplier and the horizontal axis shows the number of dimensions. The figure shows that the difference was less than a factor of 2 for the 512 dimensions with and without shared memory, but increased as the matrix size increased, reaching a maximum of approximately 2.3 times for the 4096 dimensions.

\begin{figure}[htb]
	\begin{center}
		\includegraphics[width=.75\textwidth]{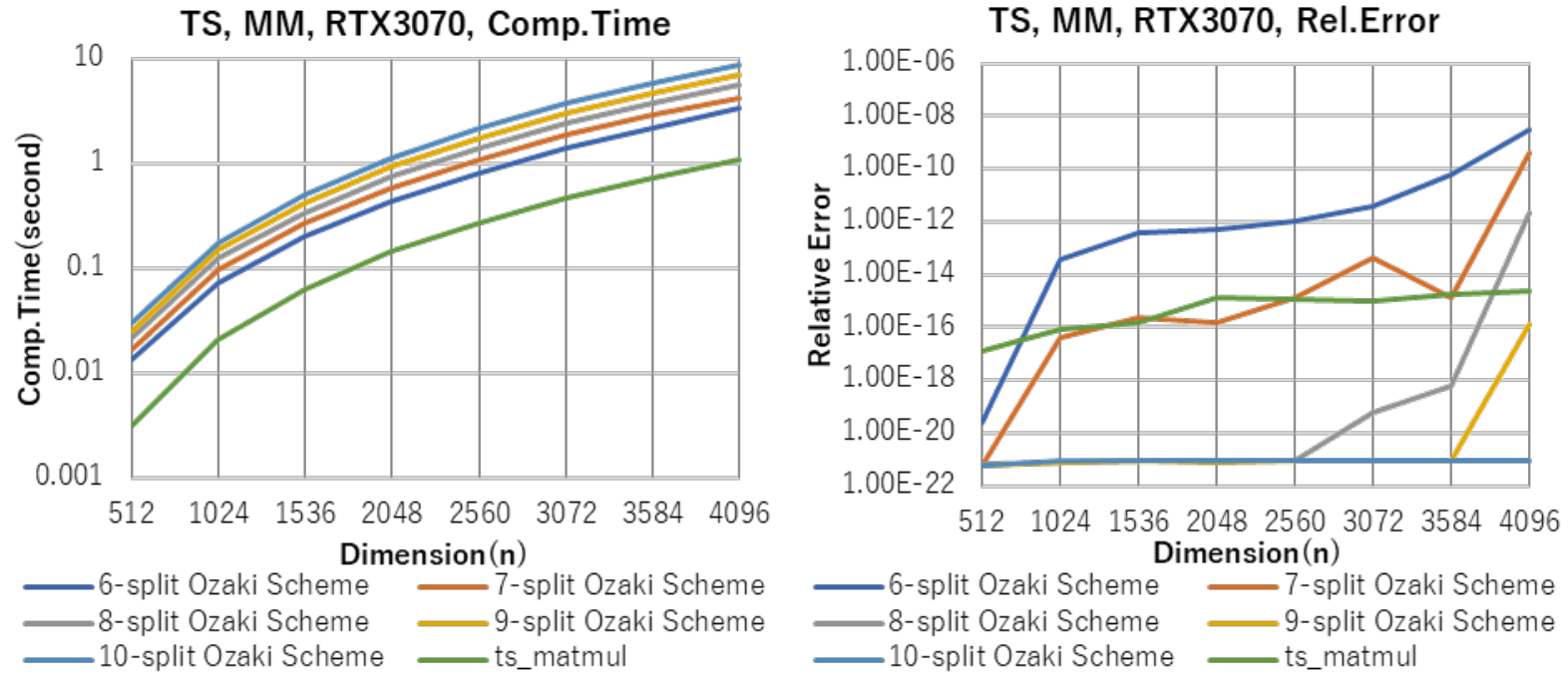}
		\caption{Computational time of TS matrix multiplication (left) and its relative error (right)}
		\label{fig:gpu_ts_time_relerr}
	\end{center}
\end{figure}

Furthermore, we evaluated the performance of Ozaki scheme on GPUs, and the results are shown in\figurename\ \ref{fig:gpu_ts_time_relerr}. For comparison, a simple matrix multiplication using shared memory is performed. The left-hand side shows the execution time and the right-hand side shows the relative error.

First, let us start with the relative error shown on the right. The vertical axis represents the relative error, and the horizontal axis represents the number of dimensions. The graph shows that the relative error does not become any smaller for Ozaki scheme with more than 10 divisions. In comparison with simple matrix multiplication, the relative error for matrix sizes 1024 to 3584 is the same for 7 divisions, and for 4096, the relative error for simple matrix multiplication falls between 8 and 9 divisions.

Next, the left figure is explained. The vertical axis represents the execution time, and the horizontal axis represents the number of dimensions. The graph shows that simple matrix multiplication is faster than the TS-type Ozaki scheme for all Ozaki schemes, suggesting that the time required by the GPU for memory allocation and access may slow down Ozaki scheme, which requires a lot of memory allocation and access. Furthermore, Ozaki scheme does not use shared memory to speed up, which may be one of the reasons for Ozaki scheme's slowness.

%
\subsection{Parallelizaion on CPU}

\begin{figure}[htb]
	\begin{center}
		\includegraphics[width=.75\textwidth]{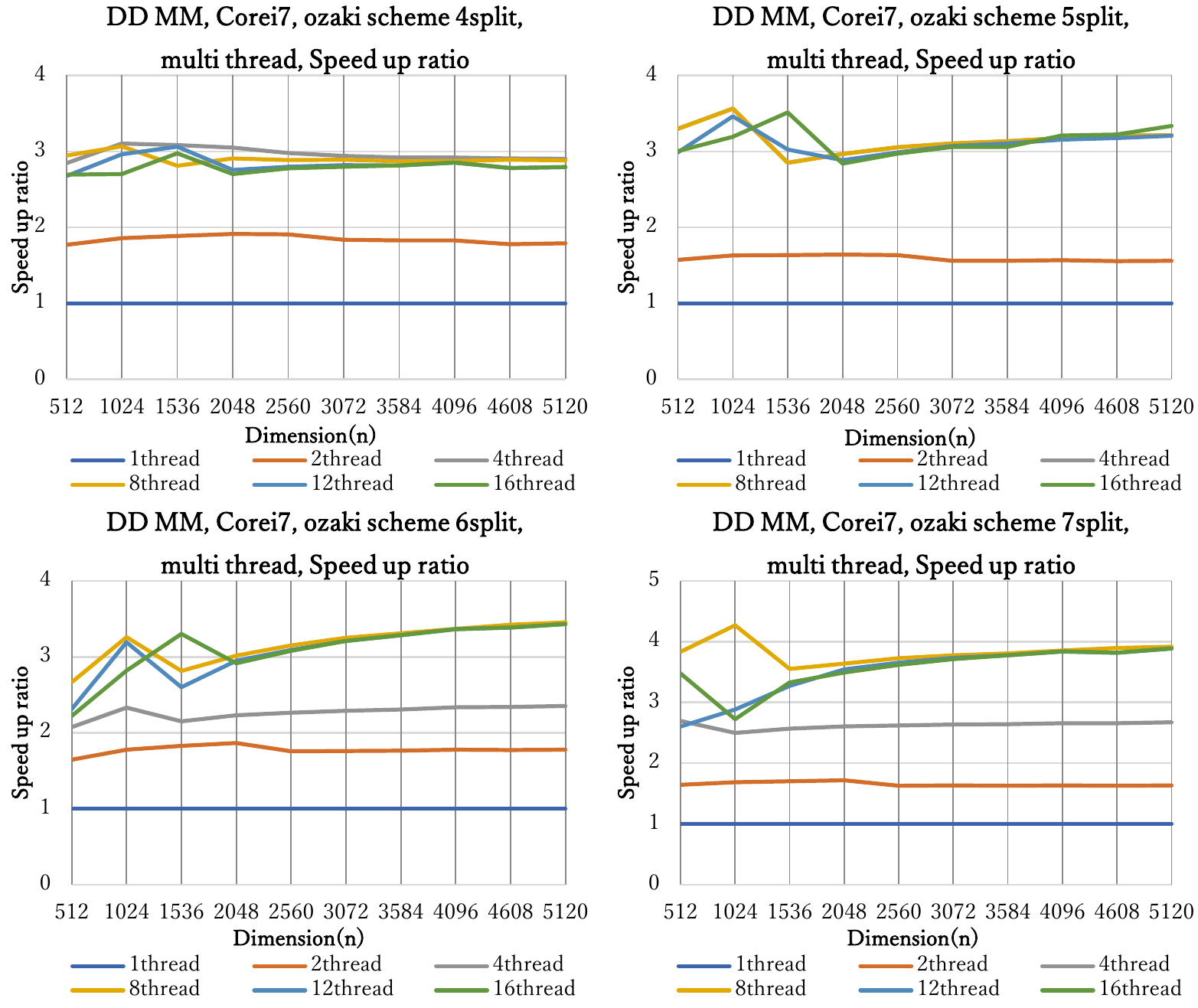}
		\caption{Computational time of parallelized DD matrix multiplication}
		\label{fig:dd_parallel}
	\end{center}
\end{figure}

\begin{figure}[htb]
	\begin{center}
		\includegraphics[width=.75\textwidth]{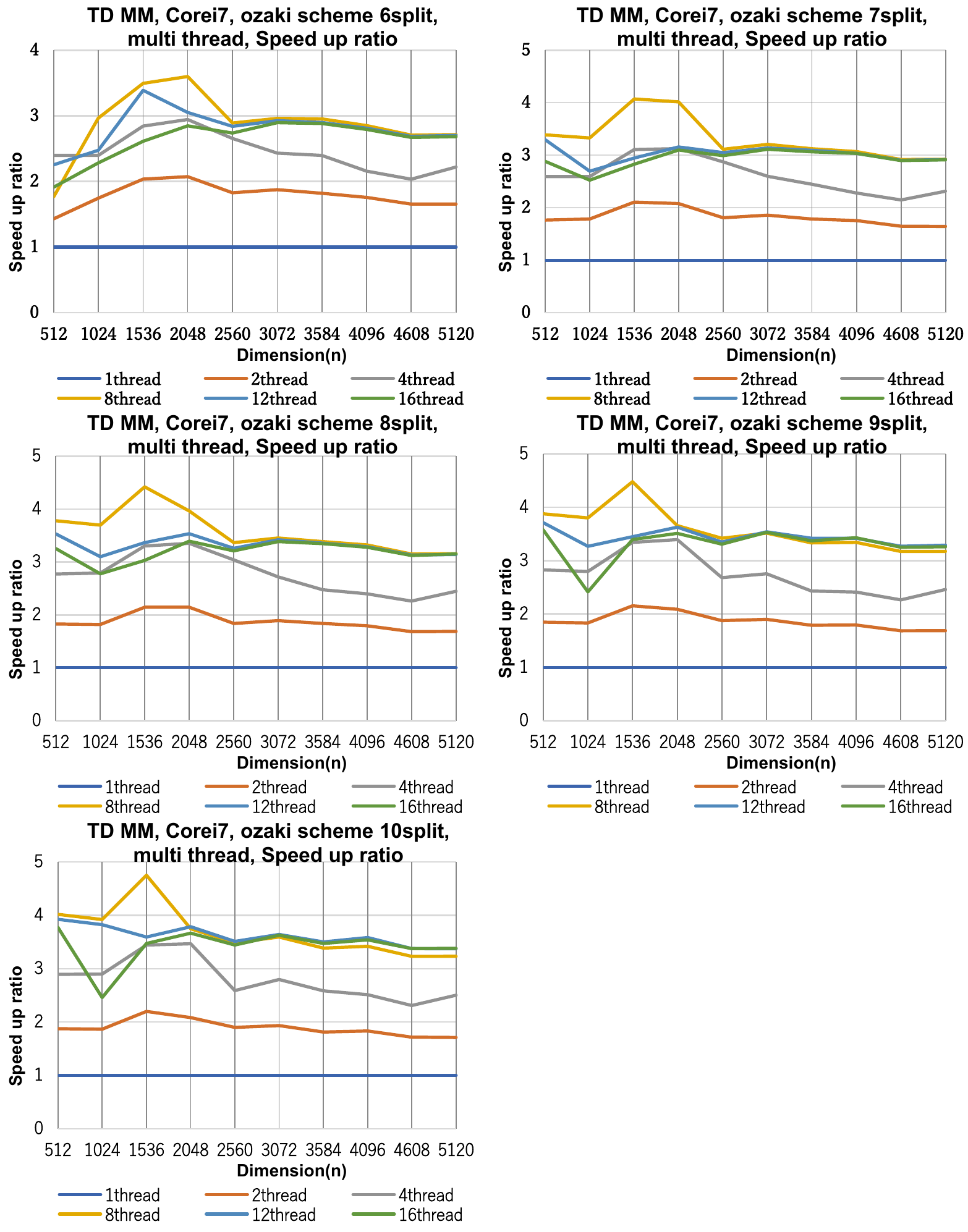}
		\caption{Computational time of parallelized TD matrix multiplication: 6 to 10 division}
		\label{fig:td_parallel}
	\end{center}
\end{figure}

\begin{figure}[htb]
	\begin{center}
		\includegraphics[width=.75\textwidth]{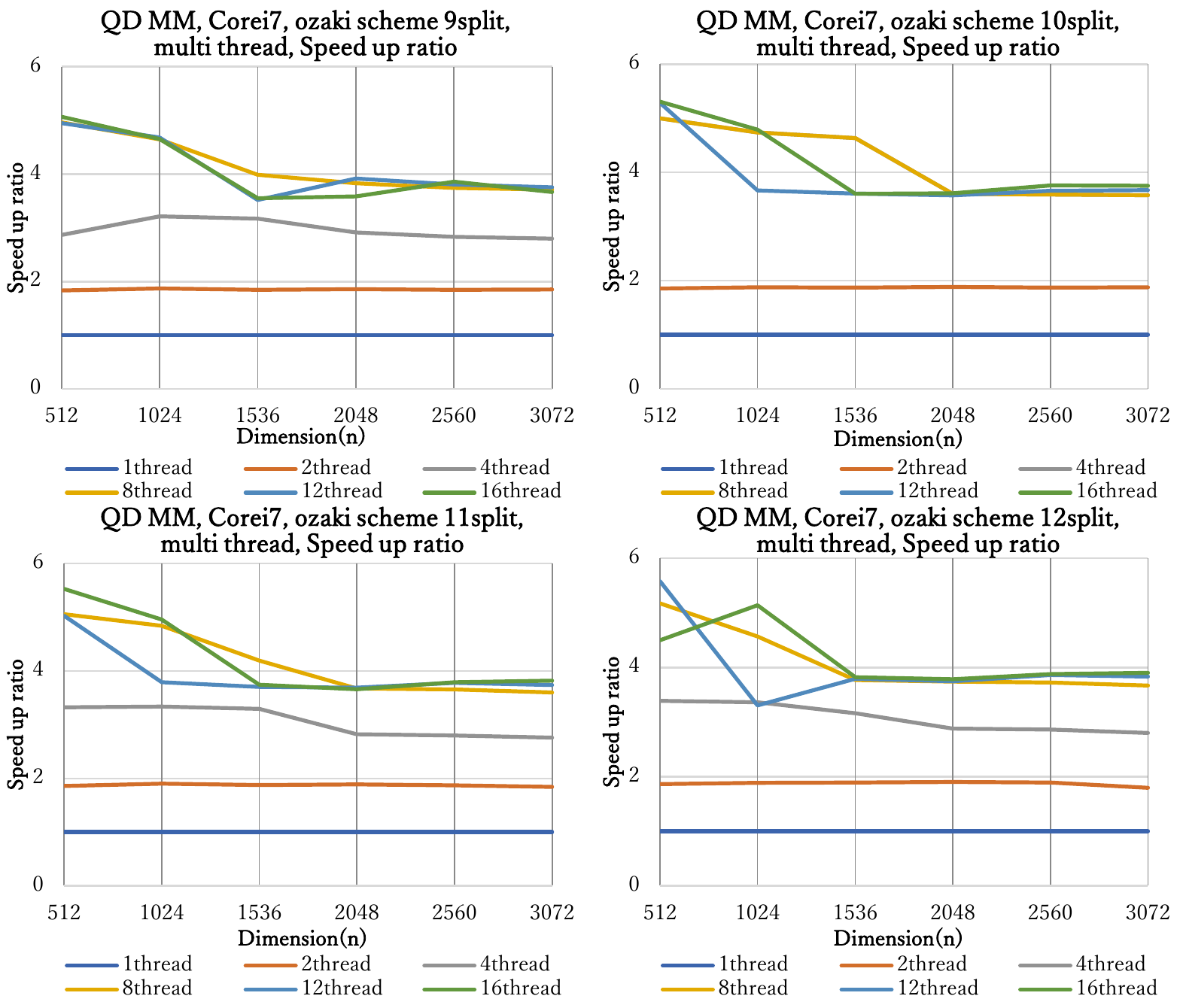}
		\caption{Computational time of parallelized QD matrix multiplication}
		\label{fig:qd_parallel}
	\end{center}
\end{figure}

The results of our parallelization with OpenMP for DD, TD, and QD matrix multiplications are shown in \figurename\ \ref{fig:dd_parallel}, \figurename\ \ref{fig:td_parallel}, \figurename\ref{fig:td_parallel}, and \figurename \ref{fig:qd_parallel}, respectively. The number of dimensions is measured at intervals of 512 dimensions from $n=512$ to $5120$. The computation time of one thread is set to 1, indicating how much the speedup is achieved when multiple threads are used. In addition, benchmarks were performed with 4 to 7 partitions of Ozaki scheme to see if the number of partitions in Ozaki scheme affects the OpenMP speedup. Parallelization was performed only for the computation of $S_A$, $S_B$, $A^{(S)}$, $B^{(S)}$, and $C_{\alpha\beta}$ in Algorithm \ref{algo:ozaki_scheme}, and for the computation loop of $C_{\alpha\beta}$, and the parallelization function provided by Intel Math Kernel was not used.

As shown in \figurename\ \ref{fig:dd_parallel}, we were able to achieve a maximum speedup of approximately 3.2 times for 8 threads with 4 divisions of Ozaki scheme and approximately 4.2 times for 8 threads with 7 divisions of Ozaki scheme in DD-type matrix multiplication. However, regardless of the number of divisions in Ozaki scheme, the speedup effect could not be confirmed for more than 8 threads. The effect of OpenMP became apparent as the number of threads in Ozaki scheme was increased, with the speedup ratio increasing by about 1 at peak time for the 7-division scheme compared to the 4-division scheme.

As show in \figurename\ \ref{fig:td_parallel}, TD-type matrix multiplication was up to about 3.5 times faster with 8 threads of Ozaki scheme with 8 divisions, and up to about 4.7 times faster with 8 threads of Ozaki scheme with 10 divisions. Regardless of the number of threads in Ozaki scheme, the speedup effect was not apparent when the number of threads was increased above 8 threads. This is because the CPU used in this study, Intel Core i7 11700, has 8 cores and 16 threads. When the number of threads in Ozaki scheme was increased, the effect of OpenMP became apparent, with a maximum of approximately 1.2 times faster with 10 threads than with 6 threads, and an overall increase of approximately 0.5 in the speedup rate.

\figurename\ \ref{fig:qd_parallel} shows that the QD-type matrix multiplication was up to 5 times faster with Ozaki scheme with 9 divisions than with Ozaki scheme with 12 divisions, and that about 5.5 times faster with 12 divisions. The effect of OpenMP was evident as the number of divisions in Ozaki scheme was increased, with the overall speedup rate increasing by about 0.5 for the 12-division scheme over the 9-division scheme.

\subsection{LU decomposition using matrix multiplication}

These results described above show that Ozaki scheme is highly effective for fixed-precision matrix multiplication, especially on CPUs. To demonstrate its potential, we present performance examples of LU decomposition in DD and TD types, which are very important in practical applications, but we would like to further optimize LU decomposition including QD and MPFR types in the future, since their application to adverse conditions is most promising.

The current LU decomposition in LAPACK is implemented using matrix multiplication, which is known to be faster than the simple LU decomposition with column-by-column computation due to the faster xGEMM. 
Therefore, if Ozaki scheme works better than Strassen matrix multiplication, it is expected to be possible to speed up the LU decomposition.

The following algorithm for LU decomposition of simultaneous linear equations can use matrix multiplication. This algorithm requires iterations. To do so, determine the width $K \in \mathbb{N}$ of the iterations.
\begin{enumerate}\small
	\item Divide $A$ into $A_{11}\in\mathbb{R}^{K\times K}$, $A_{12}\in\mathbb{R}^{K\times (n - K)}$, $A_{21}\in\mathbb{R}^{(n - K)\times K}$, and $A_{22}\in\mathbb{R}^{(n-K)\times (n-K)}$.
	\item Decompose $A_{11}$ into $L_{11} U_{11} (= A_{11})$ and then transform $A_{12}$ to $U_{12}$ and $A_{21}$ to $L_{21}$.
	\item Set $A^{(1)}_{22} := A_{22} - \underline{L_{21} U_{12}}$.
\end{enumerate}
Matrix multiplication can be used in the $L_{21} U_{12}$ part. 

\begin{figure}[htb]
	\begin{center}
		\includegraphics[width=.45\textwidth]{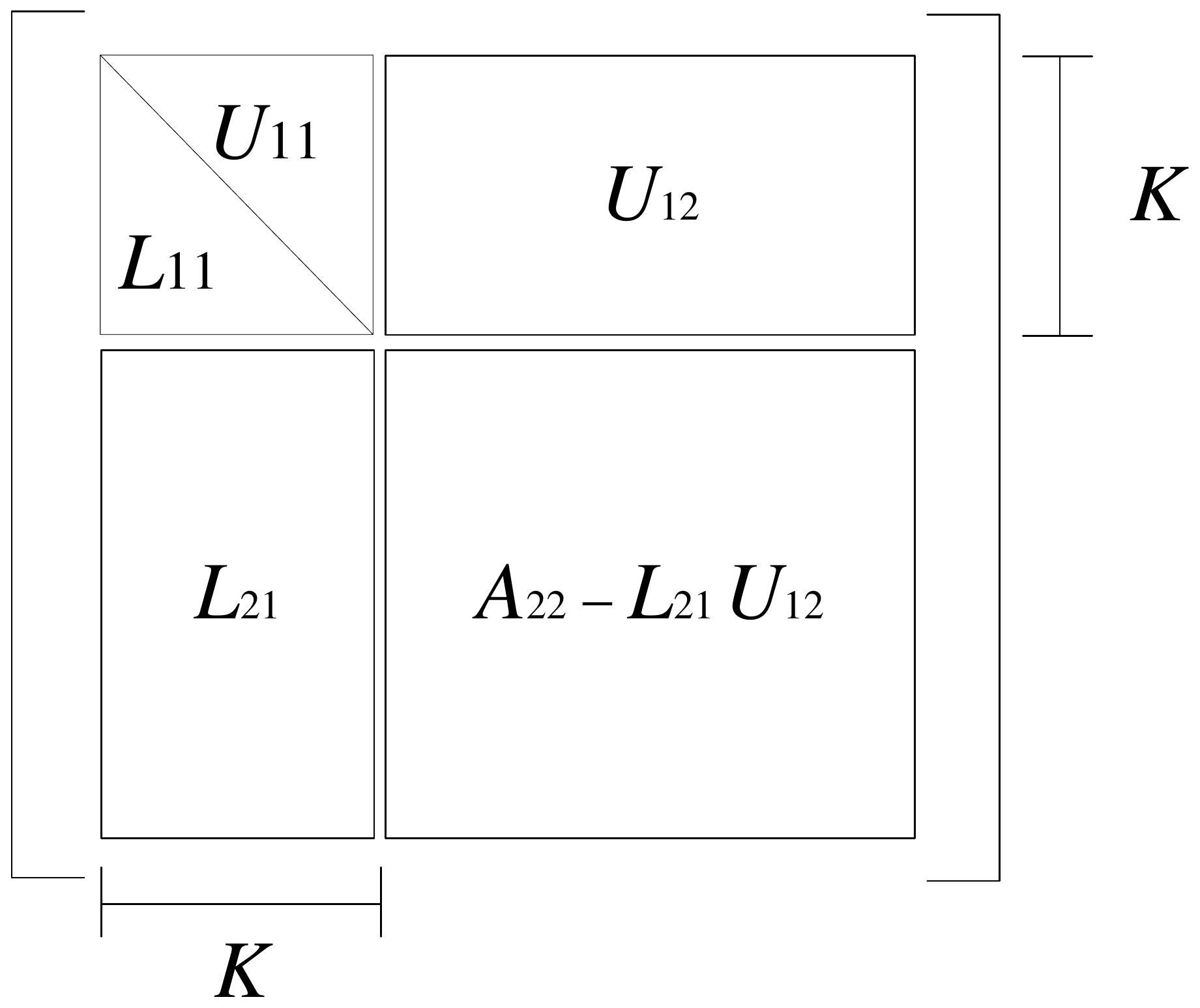}
		\caption{Algorithm of LU decomposition with xGEMM}
		\label{fig:lu_gemm}
	\end{center}
\end{figure}

The matrix multiplication was implemented using Ozaki scheme. A schematic diagram of the algorithm is shown in \figurename\ \ref{fig:lu_gemm}.
To proceed with the LU decomposition, we substitute $A= A_{22}^{(1)}$ and repeat this step as long as $n - K \geq 0$ is satisfied.

In the $n$-dimensional real coefficient linear system $A\mathbf{x} = \mathbf{b}$ used in this benchmark test, the coefficient matrix $A\mathbb{R}^{n\times n}$ uses $a_{ij}=ru$ and the constant vector $\mathbf{b}\in\ mathbb{R}^n$ uses $\mathbf{b}=[1,1/2, 1/3, ... , 1/n ]^T$. The relative error is measured as the largest relative error for each component of the numerical solution obtained by forward and backward substitution.

\subsubsection{DD-type LU decomposition}

\begin{figure}[htb]
	\begin{center}
		\includegraphics[width=.75\textwidth]{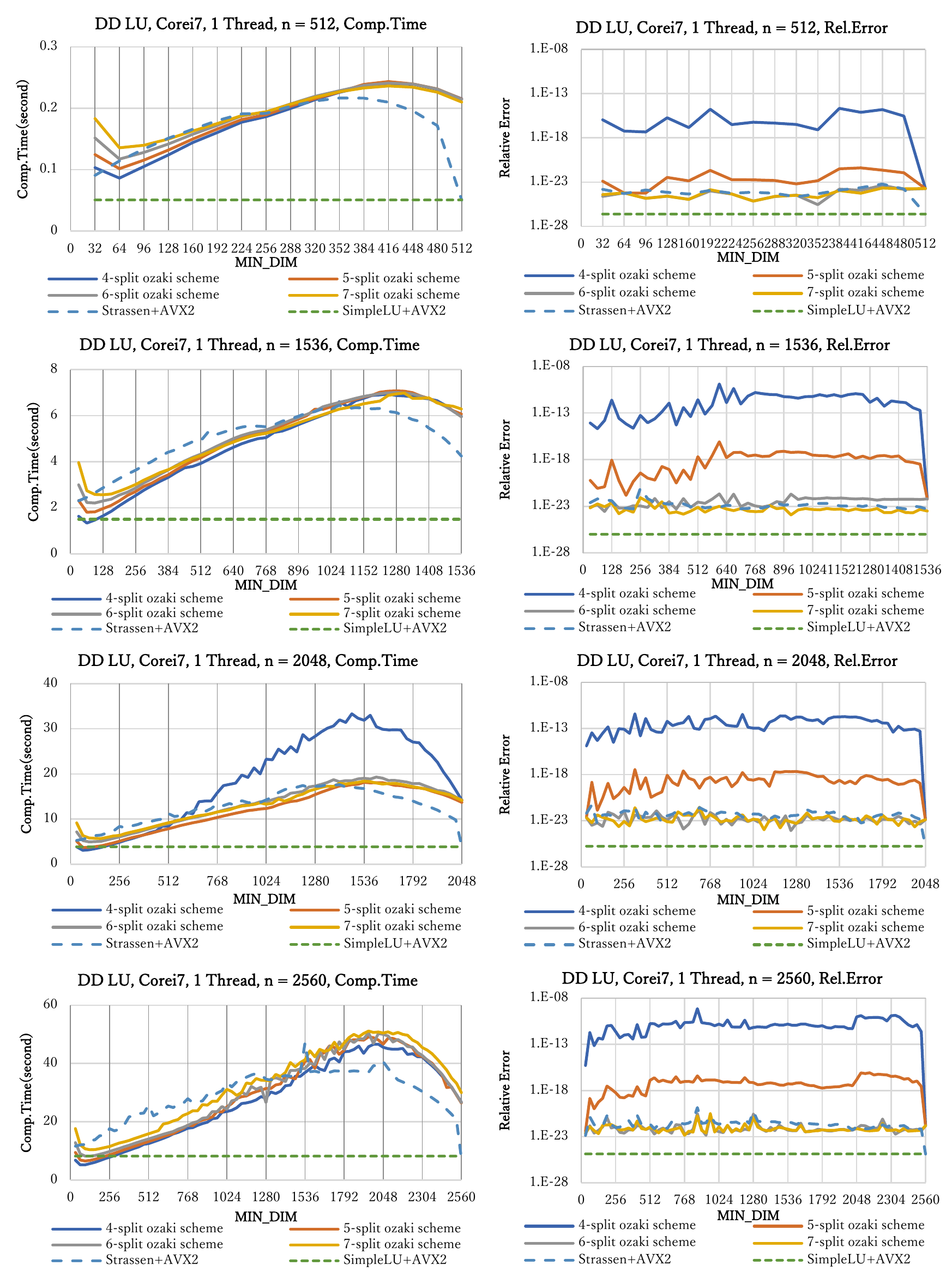}
		\caption{Computational time (left) and relative error (right) of DD LU decomposition}
		\label{fig:dd_lu}
	\end{center}
\end{figure}

The results of the performance evaluation of the LU decomposition using the DD-type Ozaki scheme are shown in\figurename\ \ref{fig:dd_lu}. LU decomposition using Strassen+AVX2 and simple LU decomposition were compared for LU decomposition using Ozaki scheme.

The left figure shows the running time and the right figure shows the relative error. The number of dimensions measured were 512, 1536, 2048, and 2560.

First, the execution time on the left side is explained. The vertical axis shows the execution time, and the horizontal axis shows the size of $K=\mathrm{MIN\_DIM}$. The size of $\mathrm{MIN\_DIM}$ is $32\leq K\leq n$, and the transition of execution time is measured by changing $K$ by 32 increments. The execution time of the simple LU decomposition is shown as a dashed line for comparison. The figure shows that the simple LU decomposition is faster than the minimum running time of Ozaki scheme 7-dimension decomposition, even for the largest matrix size measured in this study (2560 dimensions). Overall, the difference in running time between Ozaki scheme and the simple LU decomposition becomes smaller as the number of dimensions increases. Therefore, there is a possibility that the difference between Ozaki scheme and the simple LU decomposition can be reversed when the number of dimensions is larger than 2560.

Next, the relative error on the right is explained. The vertical axis represents the relative error and the horizontal axis represents the size of $\mathrm{MIN\_DIM}$. The figure shows that, as with the TD-type LU decomposition, the accuracy of Ozaki scheme is worse than that of the simple LU decomposition, by about 1 to 4 orders of magnitude. The number of divisions in Ozaki scheme with 2048 dimensions was worse than those in the other Ozaki schemes, but the reason for this was not clear.

\subsubsection{TD-type LU decomposition}

\begin{figure}[htb]
	\begin{center}
		\includegraphics[width=.75\textwidth]{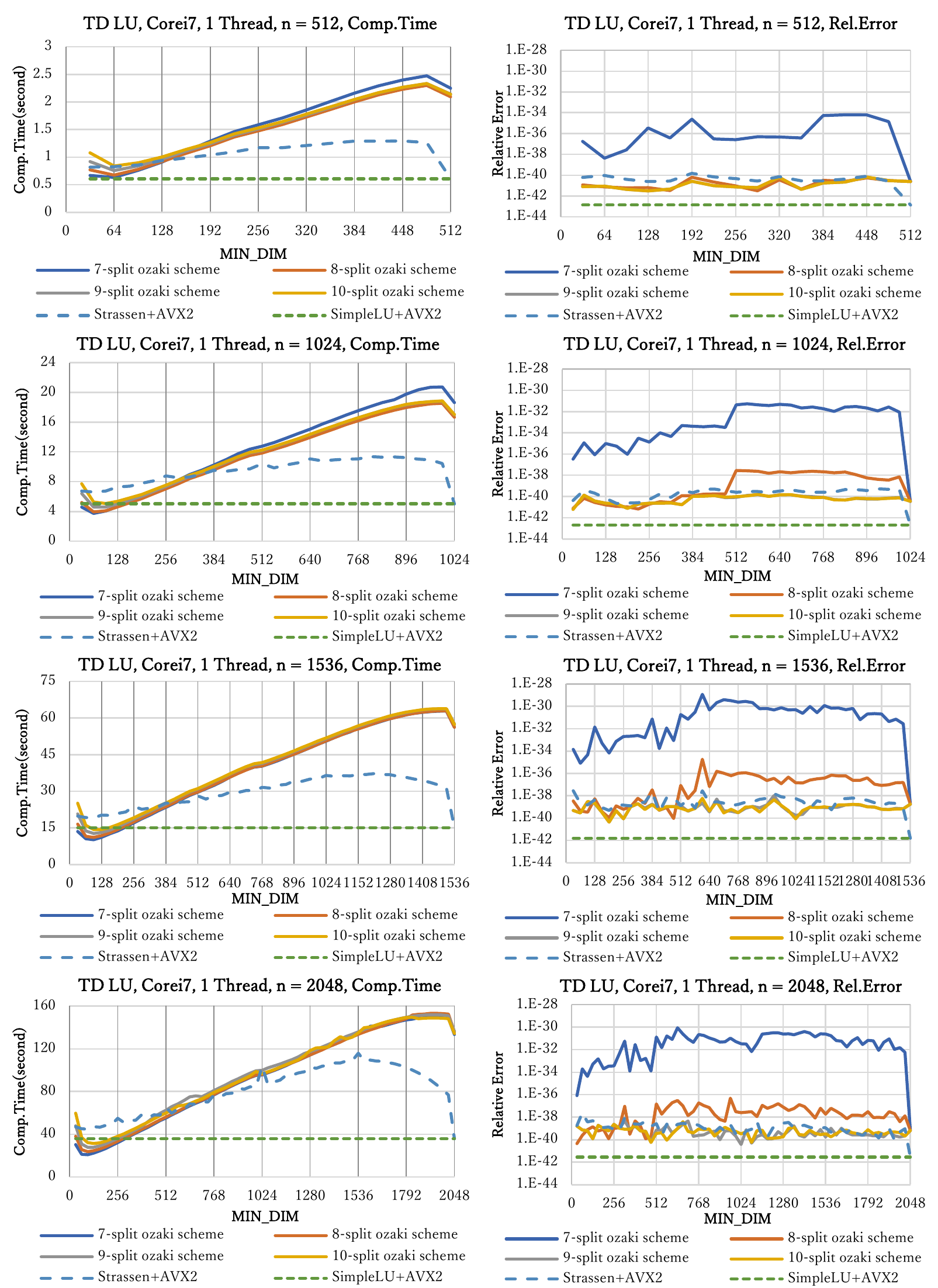}
		\caption{Computational time of parallelized TD matrix multiplication: 6 to 10 division}
		\label{fig:td_lu}
	\end{center}
\end{figure}

The results of the benchmark performance evaluation of the TD-type LU decomposition are shown in\figurename\ \ref{fig:td_lu}. As in the DD-type case, the LU decomposition using Ozaki scheme was compared to the LU decomposition using Strassen+AVX2 and the simple LU decomposition. The left figure shows the execution time and the right figure shows the relative error, measured in 512 intervals from 512 to 2048 as the number of dimensions.

First, the execution time on the left side is explained. The vertical axis of the runtime graph shows the runtime, and the horizontal axis shows the size of $K=\mathrm{MIN\_DIM}$. The size of $K$ was varied by 32 steps between $32\leq K \leq n$ and the execution time was measured. The graph shows that the simple LU decomposition is faster than the LU decomposition using Strassen+AVX2 throughout. Second, the LU decomposition using the 10-segment Ozaki scheme is faster than the simple LU decomposition from dimension 1536 onward.

Next, we explain the relative error on the right. In the relative error, the simple LU decomposition was more accurate than the LU decomposition using Strassen+AVX2 and Ozaki scheme for all $K$ in all dimensions.
In the case of $K=n$, the accuracy is better because the work is almost the same as the simple LU decomposition.
 
%
\section{Performance evaluation of arbitrary precision matrix multiplication}

These are the results of multi-component fixed-precision arithmetic in a CPU environment, showing that Ozaki scheme works well for the matrices used in this study. What about longer arbitrary precision floating-point operations? To answer this question, we introduced Ozaki scheme using DGEMM for MPFR matrix multiplication and benchmarked it against Strassen matrix multiplication by varying the number of digits $L$. The results are shown here.

\begin{figure}[htb]
	\begin{center}
		\includegraphics[width=.345\textwidth]{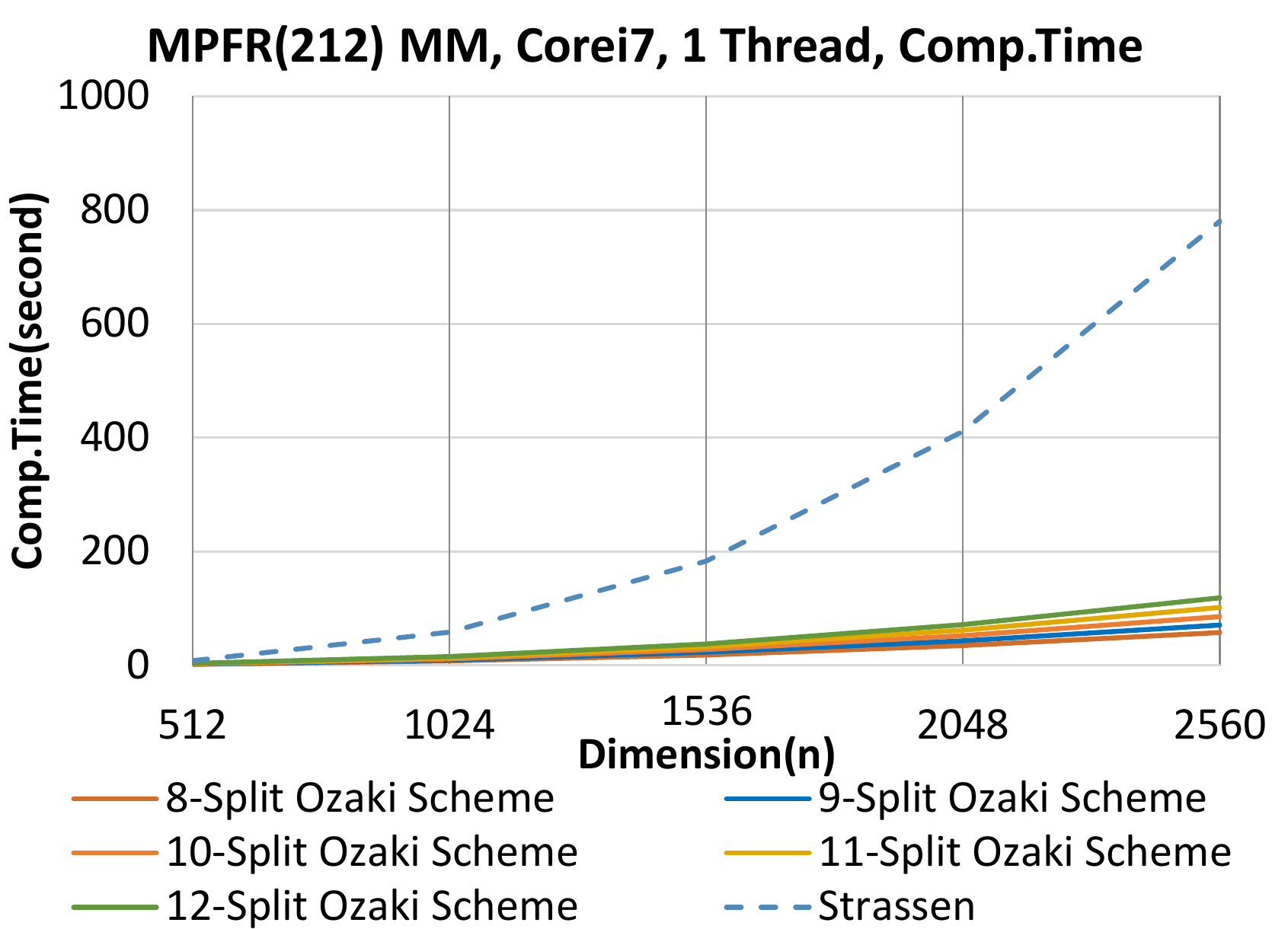}
		\includegraphics[width=.345\textwidth]{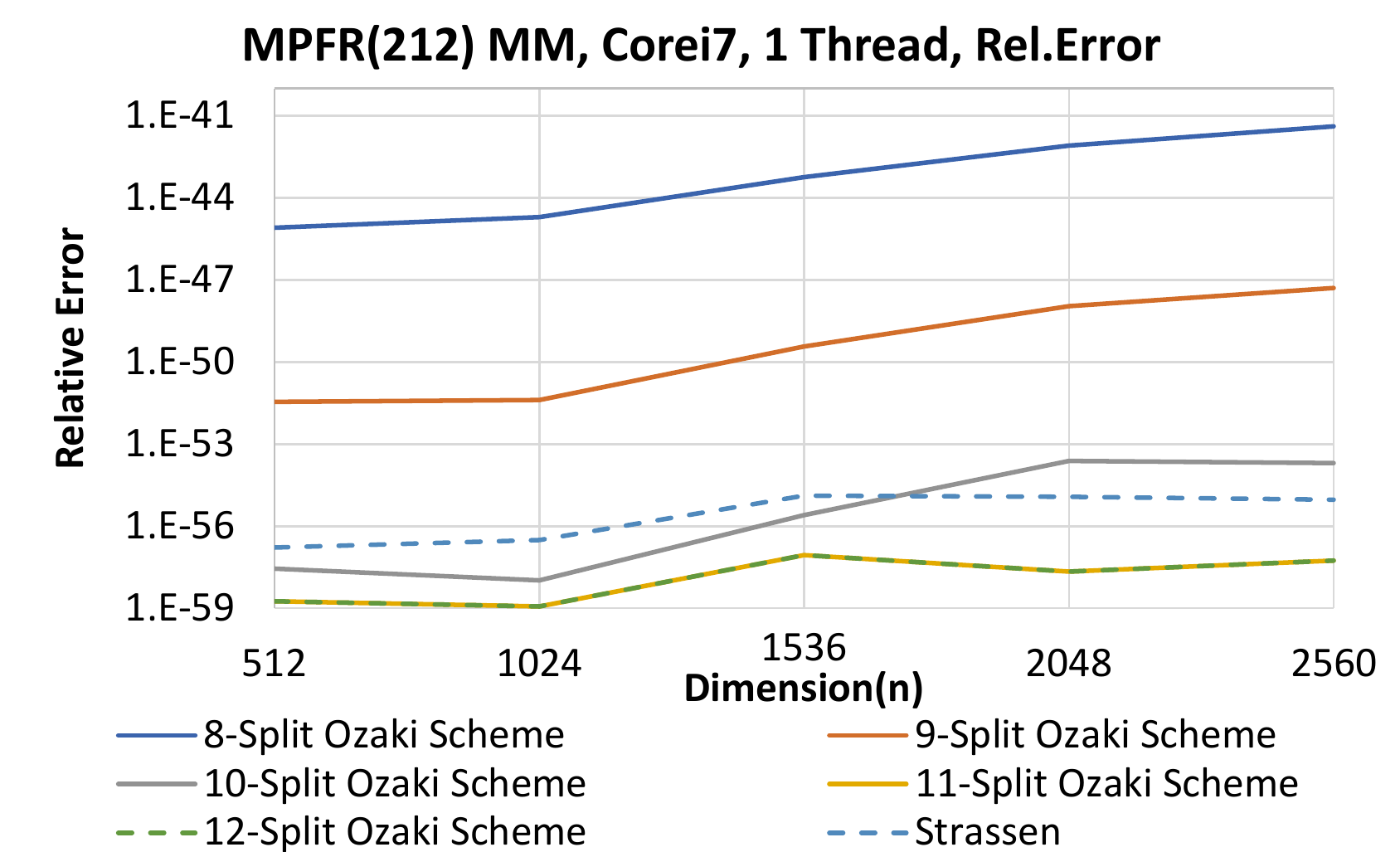}
		\includegraphics[width=.345\textwidth]{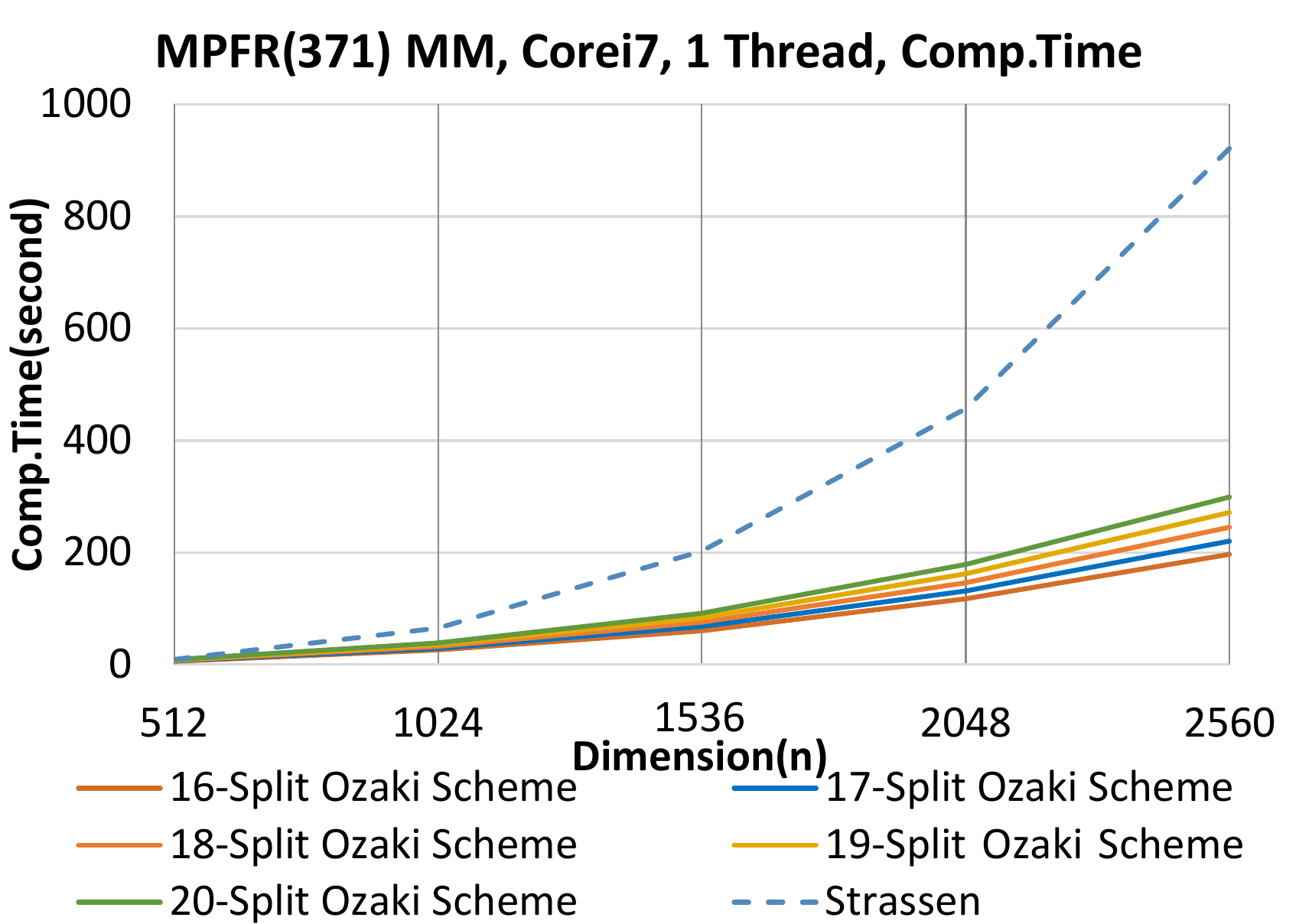}
		\includegraphics[width=.345\textwidth]{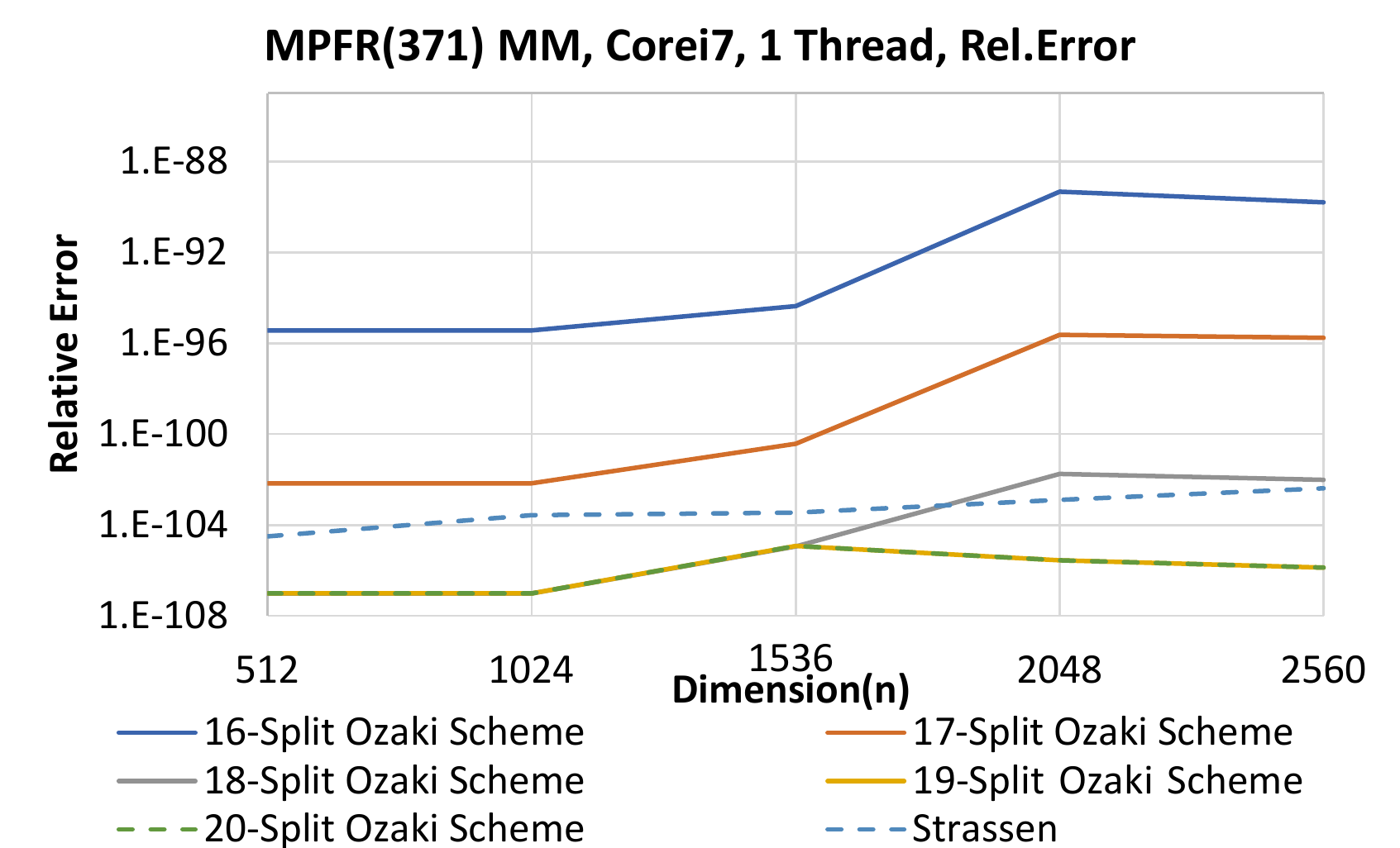}
		\includegraphics[width=.345\textwidth]{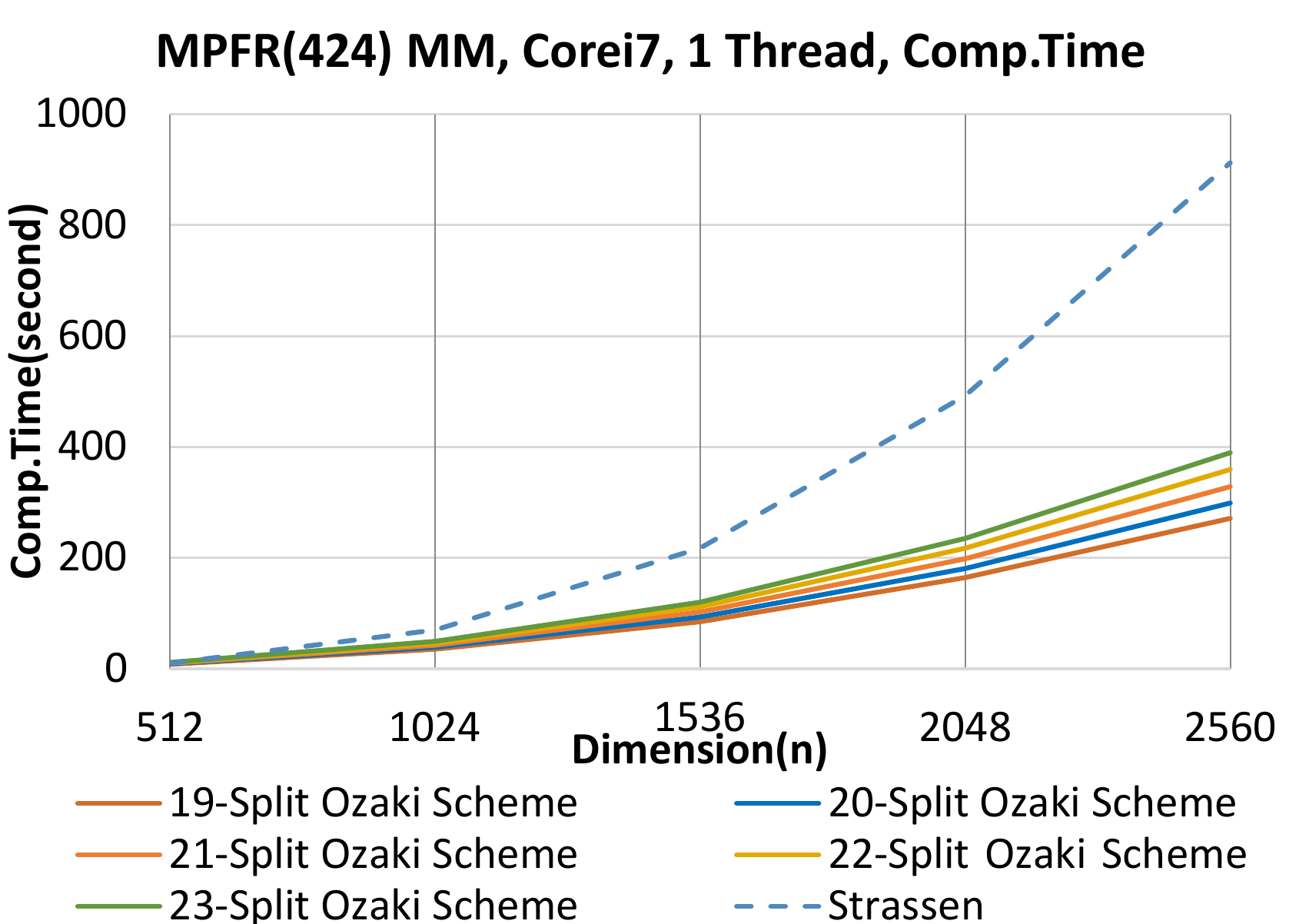}
		\includegraphics[width=.345\textwidth]{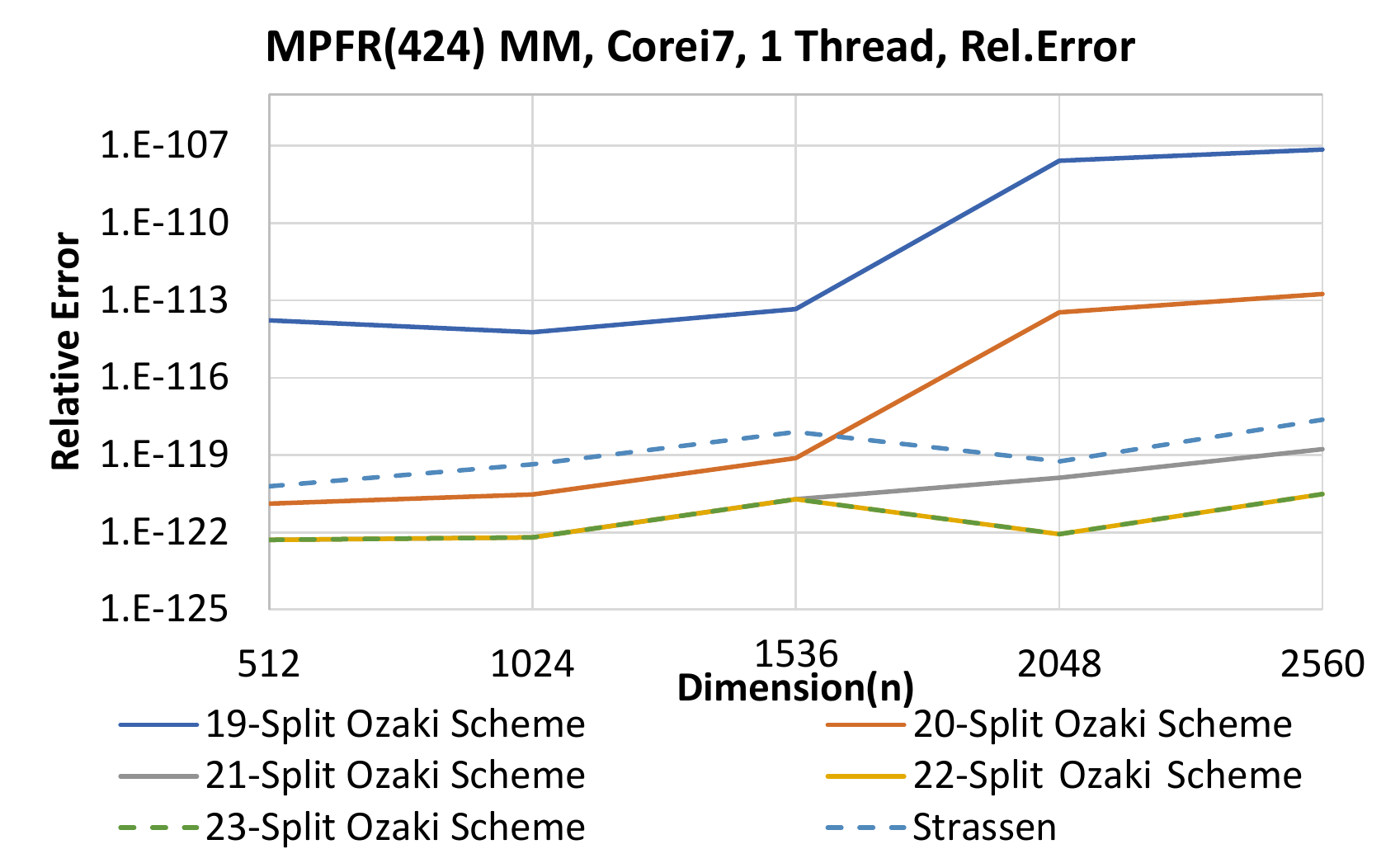}
		\includegraphics[width=.345\textwidth]{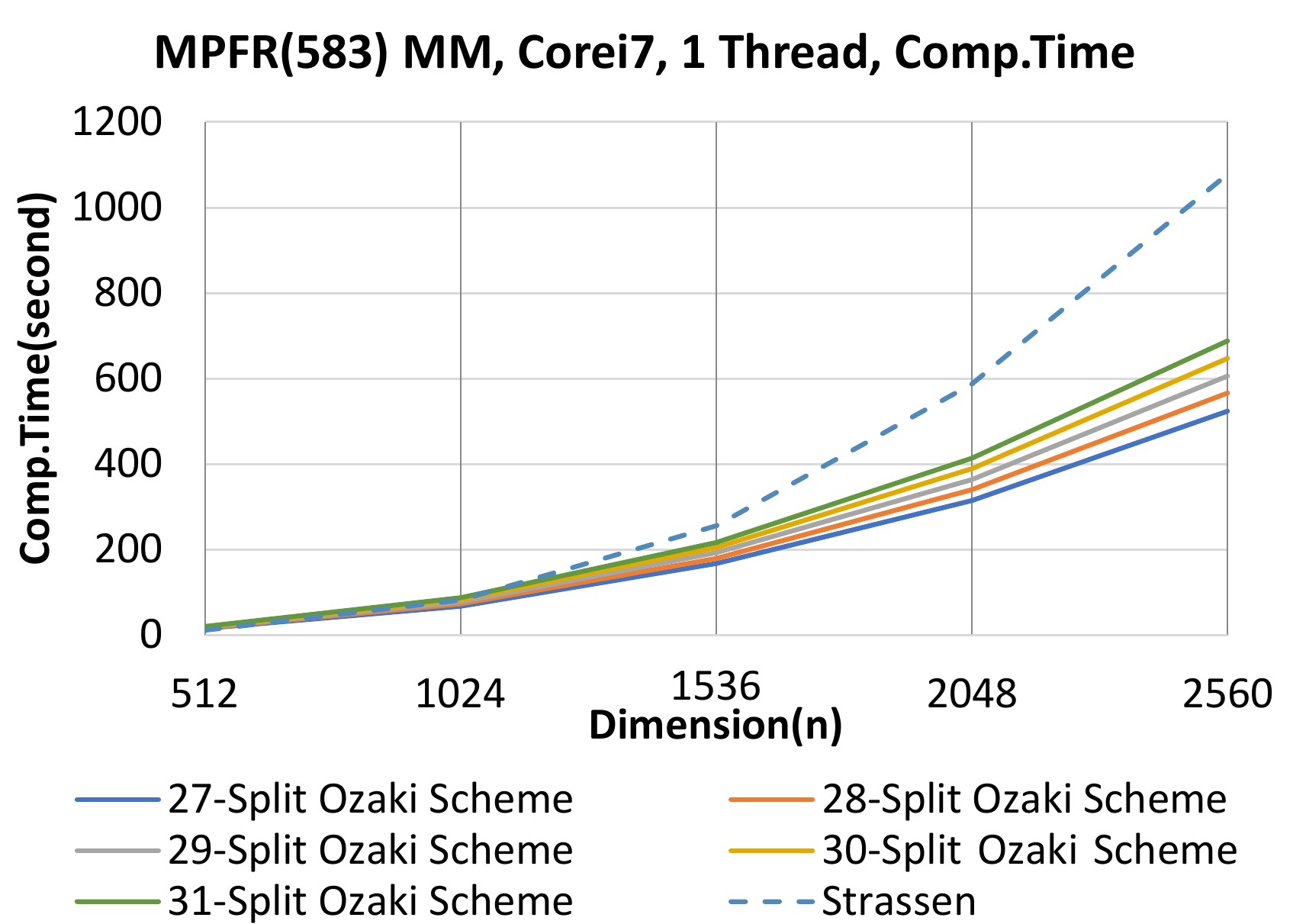}
		\includegraphics[width=.345\textwidth]{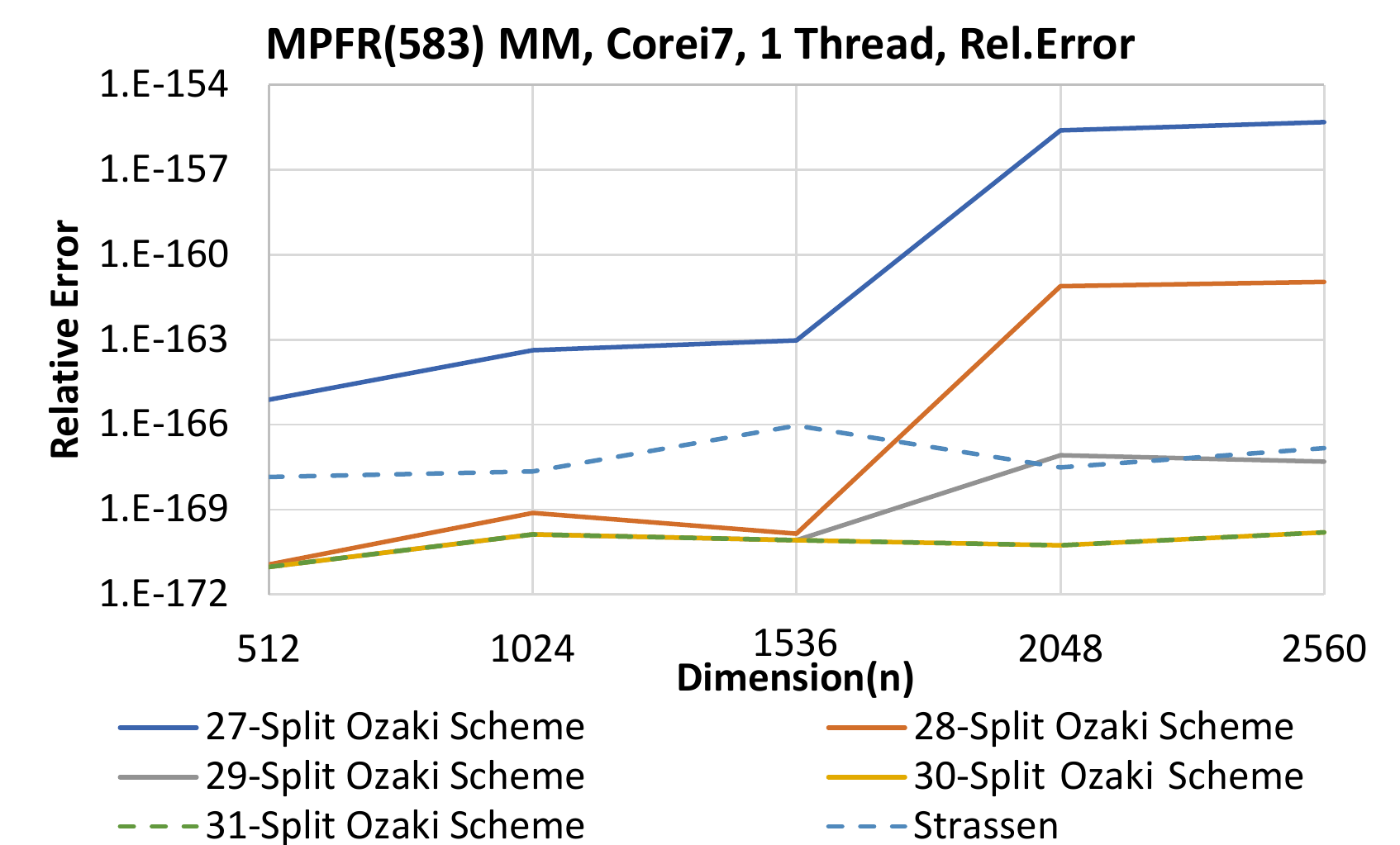}
		\includegraphics[width=.345\textwidth]{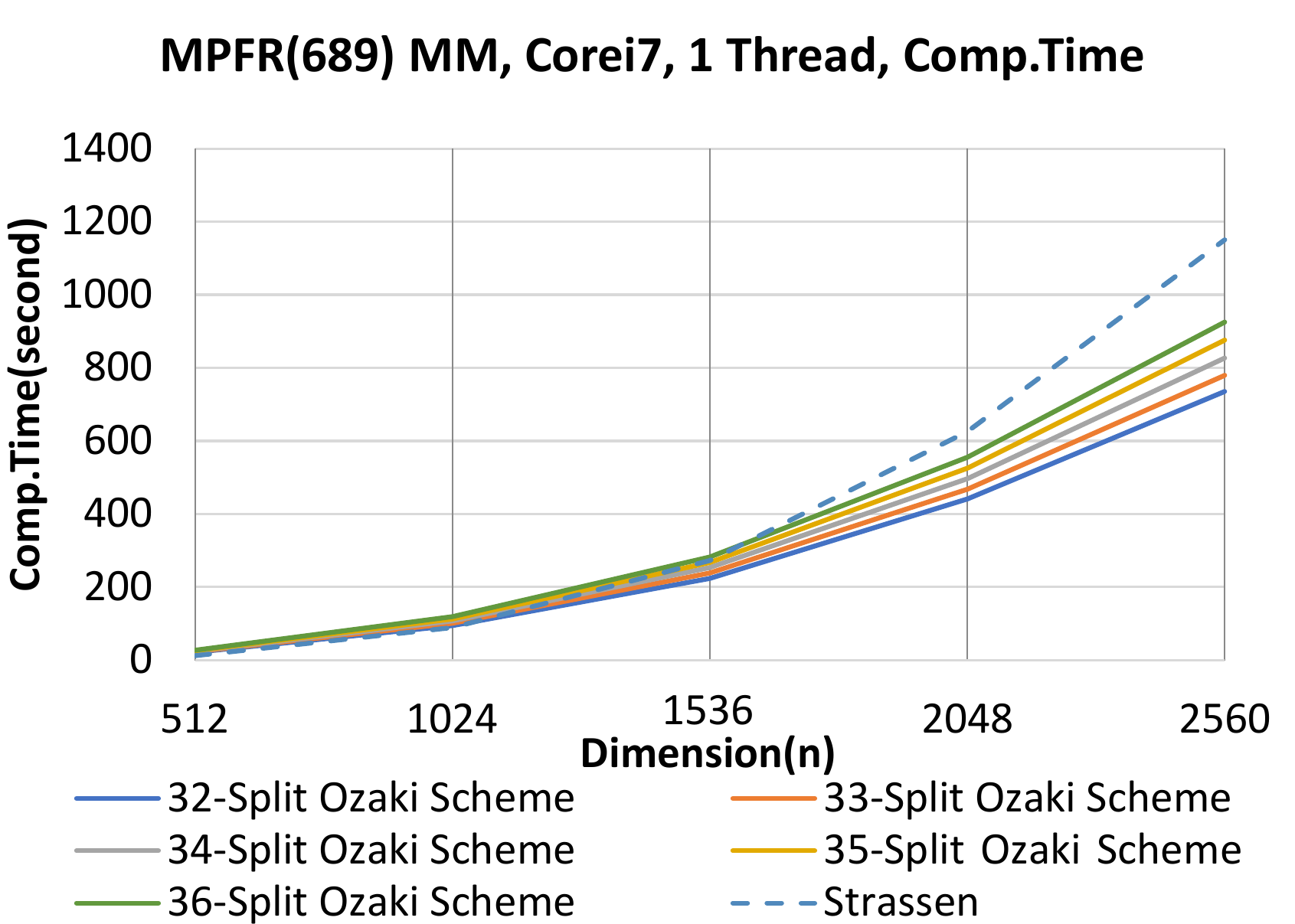}
		\includegraphics[width=.345\textwidth]{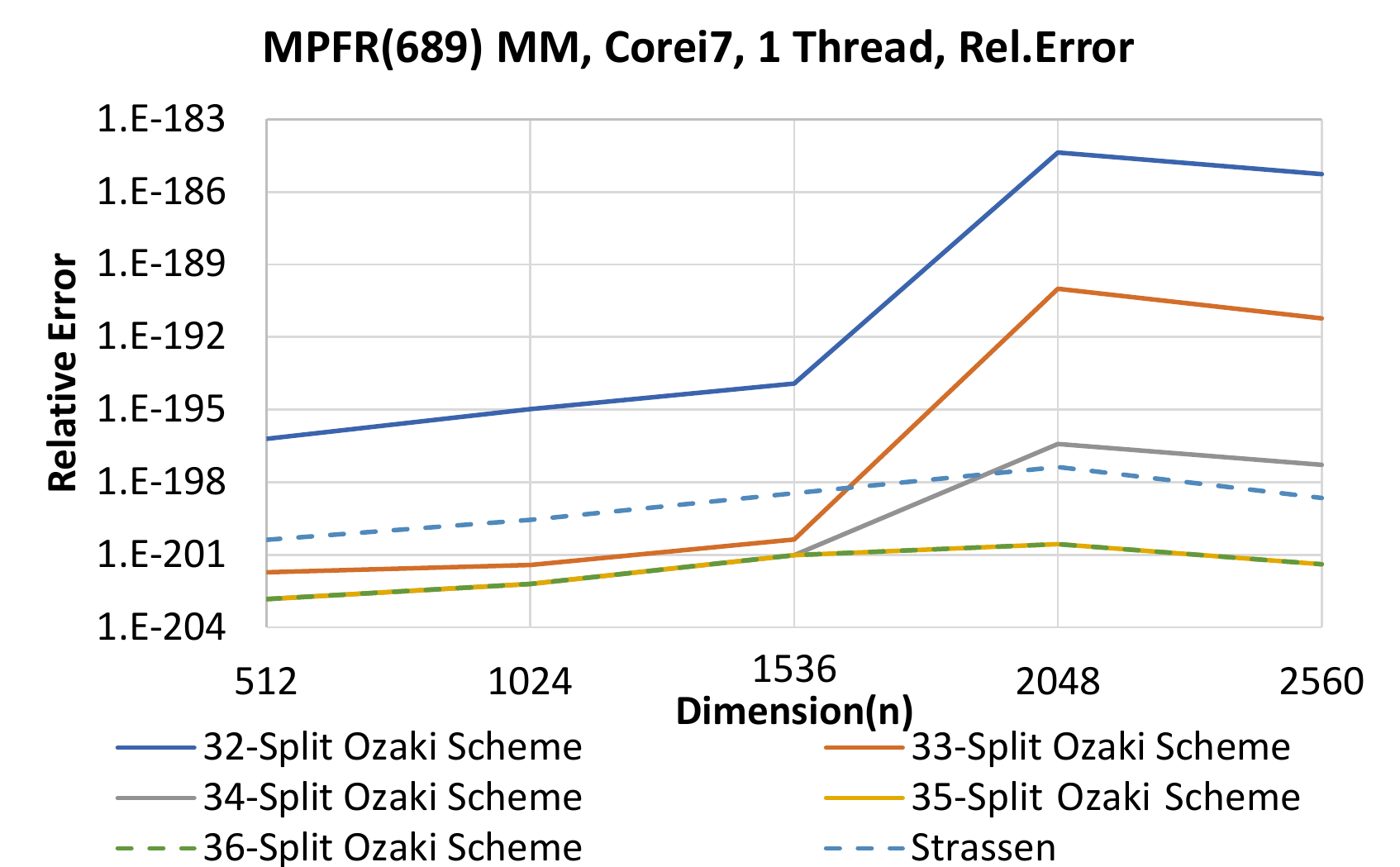}
		\includegraphics[width=.345\textwidth]{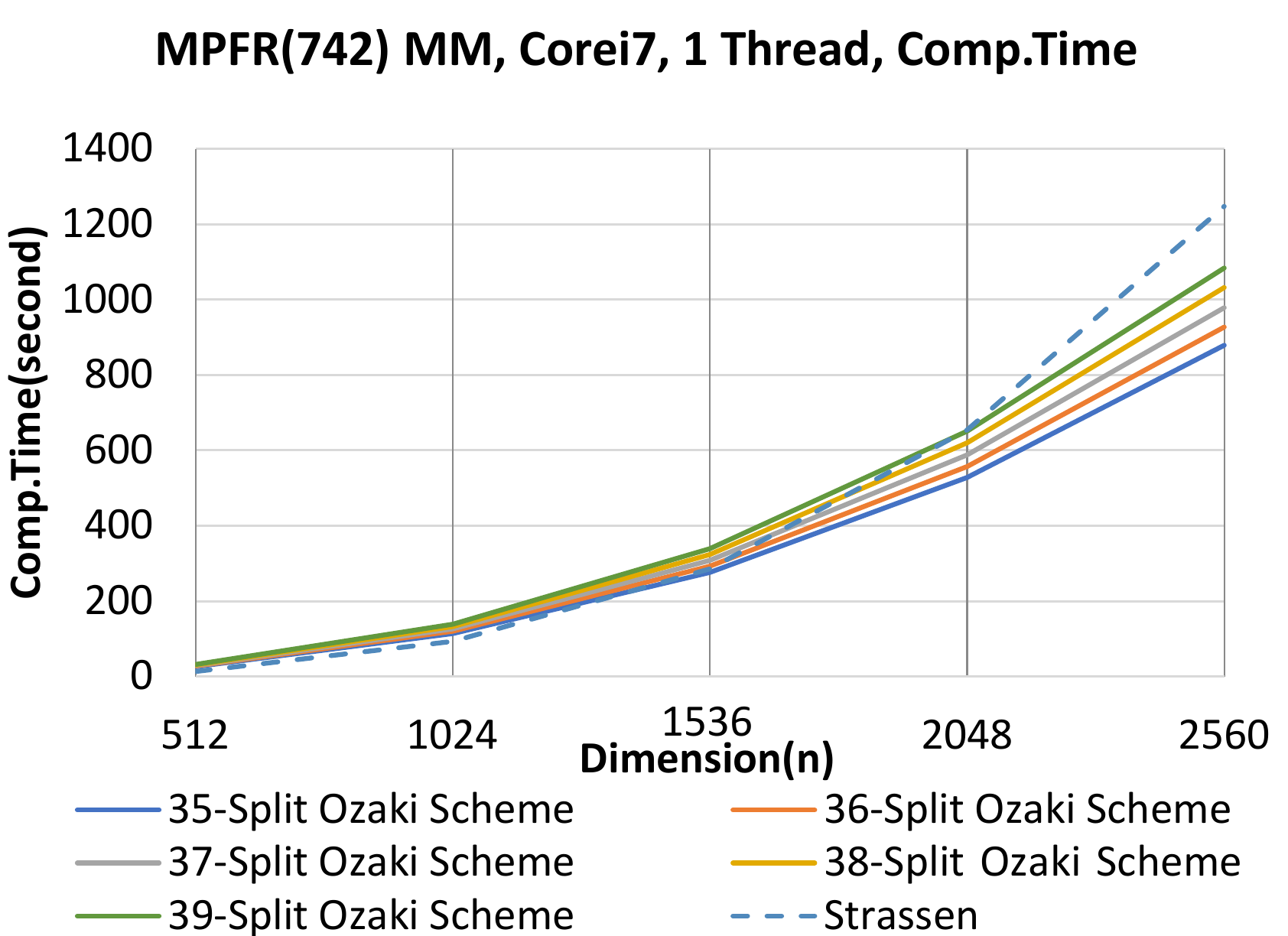}
		\includegraphics[width=.345\textwidth]{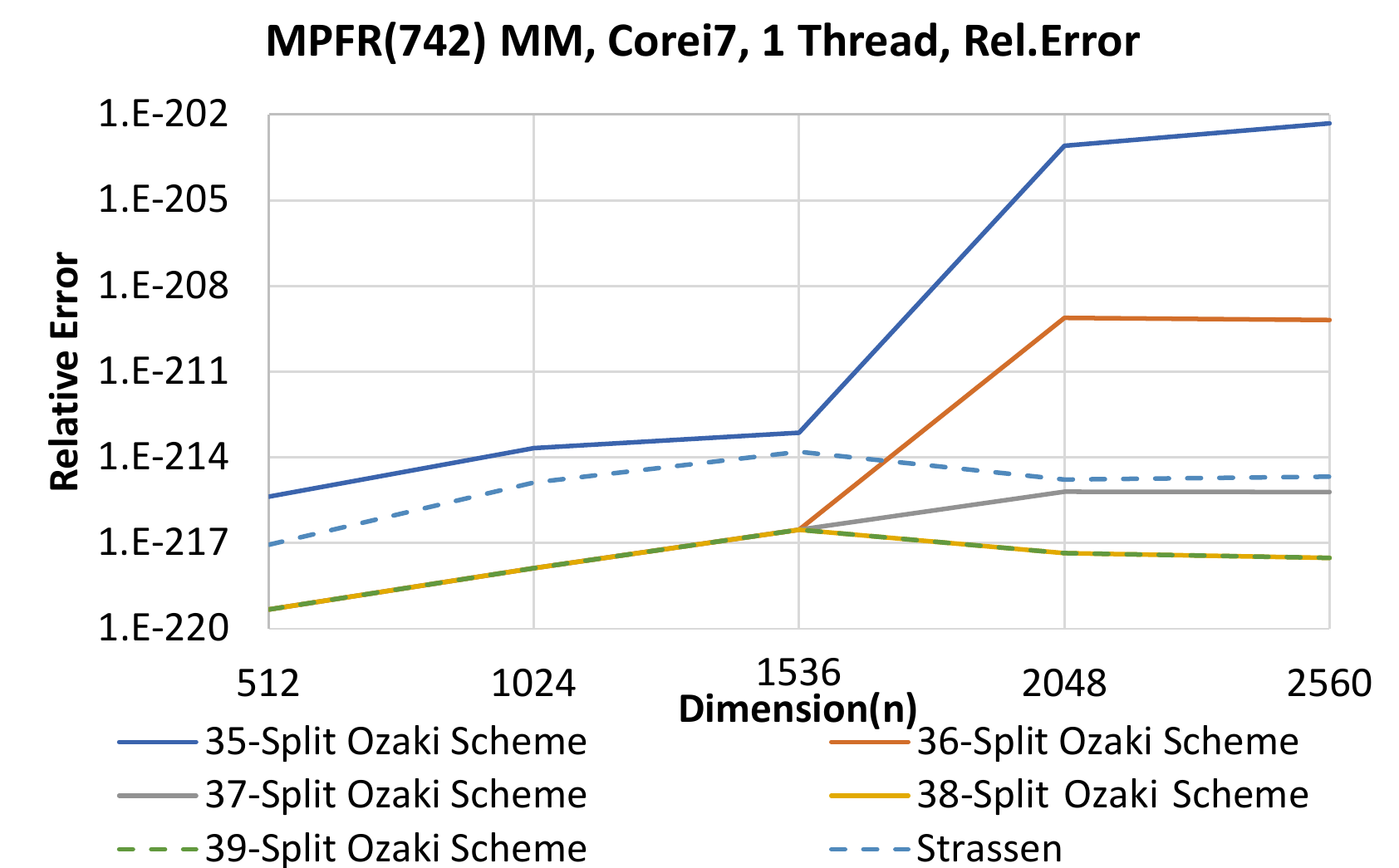}
		\caption{Computational time and relative error of MPFR 212 to 742bit matrix multiplication}
		\label{fig:mpfr_1th_time_relerr}
	\end{center}
\end{figure}

Some excerpts from the performance evaluation of matrix multiplication with the MPFR-type Ozaki scheme on a CPU are shown in \figurename\ \ref{fig:mpfr_1th_time_relerr}. The reason for using every 53 bits is to align with the mantissa part of binary64.

The library will be made available to users of multiple precision linear calculations.

From this diagram, we confirm that
\begin{enumerate}
	\item Ozaki scheme is faster than Strassen for MPFR type (424bit) up to 512 dimensions,
	\item Ozaki scheme is faster in MPFR type (583bit) up to 1024 dimensions,
	\item In 1536 dimensions, Ozaki scheme is faster up to MPFR type (689bit),
	\item In 2048 dimensions, Strassen is faster than Ozaki scheme in MPFR type (742bit).
\end{enumerate}
The results show that Strassen is faster than Ozaki's scheme in MPFR type (742bit) in 2048 dimensions. In terms of accuracy, Ozaki scheme tends to be one to four orders of magnitude more accurate than Strassen, and in terms of the number of divisions, it was found that the accuracy cannot be maintained without increasing the number of divisions by two to three for each 53-bit increase in accuracy.

The increase in the number of divisions in Ozaki scheme is one of the reasons why the increase in the number of bits shortened the execution time with Strassen. The increase in the number of divisions in Ozaki scheme is thought to have increased the number of DGEMM and MPFR additions performed internally, affecting the execution time.

%
\section{Conclusion and future work}

In this study, we implemented Ozaki scheme for DD, TD, QD, and MPFR matrix multiplication on CPUs and TS matrix multiplication on GPUs, and conducted benchmark tests. As an application of Ozaki scheme, we applied the scheme to LU decomposition for DD and TD types on CPUs and evaluated its performance.

As a result, the following results were obtained in the CPU environment:
\begin{enumerate}
	\item Ozaki scheme was approximately 1.2 times faster than Strassen+AVX2 for 7-segment and 1.8 times faster than Strassen+AVX2 for 6-segment for DD type.
	\item TD-type is about 6.9 times faster than Ozaki scheme with 9 divisions and about 8.6 times faster with 8 divisions.
	\item For the QD type, about 5.7 times faster with Ozaki scheme 11 divisions, and about 6.8 times faster with 10 divisions
\end{enumerate}
The QD-type scheme was found to be about 5.7 times faster for Ozaki scheme with 11 divisions, and about 6.8 times faster with 10 divisions.

The TD, DD, and QD matrix multiplications were parallelized using OpenMP, and the results showed that the TD, DD, and QD matrix multiplications were up to 4.7 times, 4.2 times, and 5.5 times faster, respectively, than those without OpenMP.

In the MPFR type, as the dimension increases, it becomes close to the Strassen matrix multiplication such as
\begin{enumerate}
	\item Equivalent to Ozaki scheme with maximum accuracy in 512 dimensions of the MPFR(371bit),
	\item Equivalent speed to MPFR(583) with 1024 dimensions,
	\item Equivalent speed in 1536 dimensions of MPFR(689),
	\item Equivalent speed in 2048 dimensions with MPFR(742).
\end{enumerate}
The advantage of Ozaki scheme tended to decrease as the number of bits in MPFR increased.

On the other hand, on a consumer GPU environment, we compared the TS scheme with DD and D+S schemes, and the results showed that the TS matrix multiplication scheme was approximately 9.3 times faster than DD and 11.6 times faster than D+S. However, Ozaki scheme of type TS had no advantage, and simple matrix multiplication was faster. This may be due to the overhead of function calls and the fact that Ozaki scheme does not use shared memory.

As our future works, we consider that it is necessary to determine how well LU decomposition can be adapted to QD and MPFR types, and to what accuracy LU decomposition can be implemented in QD and MPFR types, respectively. Furthermore, it is necessary to examine whether further speed-up can be achieved by changing the parallelization method, for example, by comparing the parallelization with another method, such as parallel DGEMM, other than the parallelization using OpenMP that was conducted in this study.

In parallel, we would like to implement and release a library of multiple precision basic linear computation that can easily use Ozaki scheme developed in this study to increase the convenience of multi-precision users.

%
\section*{Acknowledgement}
This research was supported by Grant-in-Aid for Scientific Research 20K11843. It was also partially supported by Shizuoka University of Science and Technology Proposal Research Fund. Personal communication at HPC research meetings helped to improve the content. We thank all parties involved.


\end{document}